\documentclass[]{interact}

\usepackage{graphicx}
\usepackage[utf8]{inputenc}
\usepackage{mathtools}
\usepackage{multirow}
\usepackage{graphicx}
\usepackage{subfigure}
\usepackage{hyperref}
\usepackage{apacite}
\usepackage{natbib}

\begin{document}

\title{A bi-integrated model for coupling lot-sizing and cutting-stock problems}

\author{
\name{Amanda O. C. Ayres\textsuperscript{*}\thanks{\textsuperscript{*}Corresponding author: amanda@dca.fee.unicamp.br}, Betania S. C. Campello, Washington A. Oliveira \and Carla T. L. S. Ghidini}
\affil{School of Applied Sciences, University of Campinas, R. Pedro Zaccaria, 1300, Limeira, 13484-350, SP, Brazil}
}

\maketitle

\begin{abstract}
In this paper, a framework that addresses the core of the papermaking process is proposed, starting from the production of jumbos and ending with the paper sheets used in daily life. The first phase of the process is modelled according to a lot-sizing problem, where the quantities of jumbos are determined in order to meet the demand of the entire chain. The second phase follows a one-dimensional cutting-stock formulation, where these jumbos are cut into smaller reels of predetermined lengths. Some of these are intended to fulfil a portfolio of orders, while others are used as raw material for the third phase of the process, when the reels are cut into sheets with specific dimensions and demands, following a two-dimensional cutting-stock problem. The model is called the Bi-Integrated Model, since it is composed of two integrated models. The heuristic method developed uses the Simplex Method with column generation for the two cutting-stock phases and applies the Relax-and-Fix technique to obtain the rounded-integer solution. Computational experiments comparing the solutions of the Bi-Integrated Model to other strategies of modelling the production process indicate average cost gains reaching 26.63\%. Additional analyses of the model behaviour under several situations resulted in remarkable findings.
\end{abstract}


\section{Introduction}

\begin{figure}[!h]
	\centering
	\includegraphics[width=\textwidth]{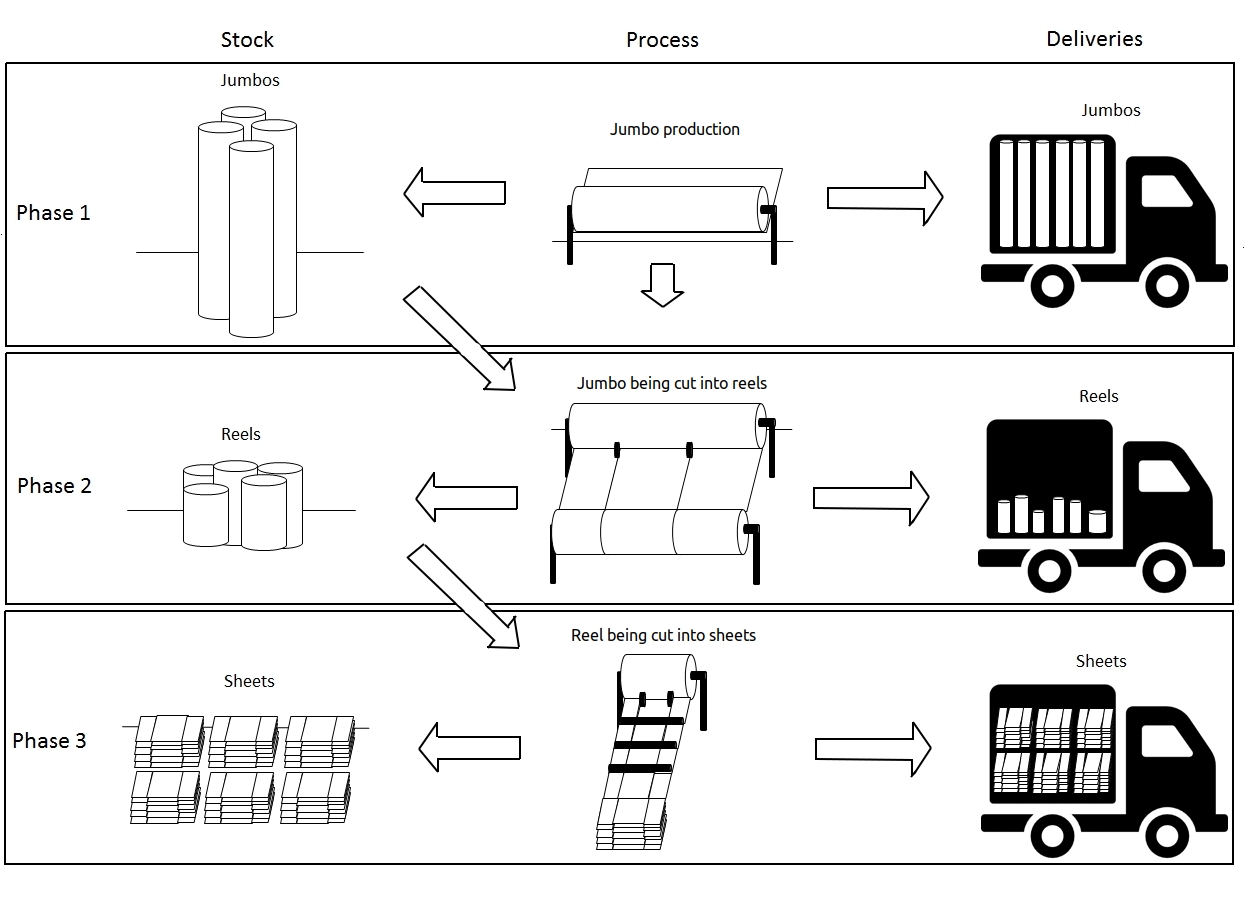}
	\caption{The production of jumbos, reels and sheets.}
	\label{figure:process}
\end{figure}

Paper production is an intricate and large-scale process consisting of three main phases, as illustrated in Figure \ref{figure:process}. In Phase 1, cellulose sheets are rolled around an axis, resulting in large cylinders called jumbos. In Phase 2, machines known as rewinders unwind these jumbos, make longitudinal cuts and rewind them continuously. The smaller barrels resulting from this process are called reels. In Phase 3, the reels are processed by machines called cutters that unwind them while making longitudinal and transverse cuts, producing the paper sheets used in daily life. All three products---jumbos, reels and sheets---can either be delivered to the customer, stored as stock or served as raw material to the next phase of the process, if any. They all share a critical feature called grammage, which is the mass of paper per unit of area and defines the product type. Once the grammage is set for the jumbo at the beginning of the process, the entire chain will maintain the same specification. The quantity of machines in each phase varies according to the size of the factory, but it generally consists of some units.

There are basically three types of costs involved in each stage of the process---production, machine setup and storage costs. Whenever there is a change in the product specifications (grammage or dimensions of the final items), the three types of machines on the production line must be reset, to change the positions of the knives responsible for cutting in each phase. The cost of this wasted time is accounted for as the machine setup cost and can be estimated as the profits of the products that could be manufactured during this period.

Several works address this production process. \citet{Poltroniere2008a}, \citet{Poltroniere2016}, \citet{Leao2017250}, \citet{Malik} and \citet{Campello2017} integrated Phases 1 and 2 using an approach that deals with lot-sizing and cutting-stock problems. In \citet{Poltroniere2008a} and \citet{Poltroniere2016}, machine capacity and setup constraints were considered in the formulations. This model interconnected the phases by integrating the constraints, whereby the quantity of jumbos manufactured was linked to the quantity needed to meet the demand for reels. \citet{Leao2017250} presented three mathematical formulations---the first oriented to items, the second to cutting patterns and the third to machine decomposition. \citet{Malik} proposed a model to jointly optimize the changeover, inventory holding cost and trim loss, also discussing the relationship between the ensuing cycle service levels and the total joint costs. \citet{Keskinocak2002} also integrated Phases 1 and 2, but focused on allocation and sequencing problems. In addition to jumbo and reel production programming, reel distribution was also considered using a multi-criteria approach. \citet{Campello2017} used a multi-objective approach to address the two problems being integrated.

Dealing with Phase 2, \citet{Schrage2002} modelled the allocation of orders to reels cutters using the classical assignment problem. \citet{Gomory1961}, \citet{Zak2002} and \citet{Cristina2016} used the cutting-stock problem in their theoretical work. In the well-known work of \citet{Gomory1961}, the problem of one-dimensional jumbos being cut into reels was solved using the Simplex Method with column generation. \cite{respicio2002integrating} extended the idea in \citet{Gomory1961} to a case where the production capacity is not enough to meet the demand, resulting in long delivery times. \citet{Zak2002} generalised this approach for multi-stage, one-dimensional cutting. \citet{Cristina2016} assumed jumbo production to be an input parameter, and stock balance constraints required planning for jumbo stock in each period, thus allowing reel production to be anticipated.

Integrating Phases 2 and 3, \citet{Chauhan2008} used a classic assignment problem. The proposed decision model determined the optimal quantity of reels to be stored and the corresponding number of sheets according to specific demands, resulting in a smaller loss of material. \citet{Correia2004a} and \citet{Kallrath2014a} also integrated the two phases using a cutting-stock approach. \citet{Correia2004a} enumerated the combinations of reel widths and determined the quantities to be produced and the cutting patterns to be used by employing a mathematical model for the 1.5-dimensional cutting-stock problem. \citet{Kallrath2014a} presented several approaches to solve one- and two-dimensional cutting-stock problems, including a new column generation heuristic method of treating different jumbo widths, with limited jumbo stocking and sub-production allowed but to be avoided. \citet{Krichagina} considered lot-sizing and cutting-stock problems to produce sheets of several sizes to minimise waste costs, machine setup and sheet stock.

\citet{Vecchietti2008}, \citet{Gomory1965} and \citet{Kim2014} addressed Phase 3 only. \citet{Vecchietti2008} considered allocation and cutting-stock problems to optimise the production of corrugated paper boxes. \citet{Gomory1965} generalised the work of \citet{Gomory1961}, addressing the staged multi-dimensional cutting-stock problem. \citet{Kim2014} worked with two-stage two-dimensional cutting-stock problems, considering several factory operation constraints.

\citet{melega2018classification} presented a survey about integrated lot-sizing and cutting-stock problems in several types of industries, such as textiles, wood processing, aluminium, copper and furniture. They also proposed a general model to encompass all the cited models and classify them according to several criteria.

The present paper proposes a novel formulation that integrates the three phases of paper production in a unique and general model that can be fitted to any factory. The first phase is modelled according to a lot-sizing problem, the second phase is modelled according to a uni-dimensional cutting-stock problem and the third phase is modelled according to a bi-dimensional cutting-stock problem. It is important to note that, although this work has a similar proposal to the paper by \citet{melega2018classification}, in the sense of integrating the three stages of production into a generalist model, \citet{melega2018classification} presented a model where Phases 1 and 3 deal with lot-sizing problems, while Phase 2 deals with a cutting-stock problem. In our work, Phase 1 deals with a lot-sizing problem, while Phases 2 and 3 deal with cutting-stock problems. This feature is specific of the paper industry and, to the best of our knowledge, this is the first time it has been addressed in the literature.

We compare the results of this new integrated model with those obtained by employing other strategies to approach the overall process. Furthermore, we explore the model's behaviour in several scenarios and present some striking findings. It is the first time that a model has dealt with the three phases of paper production together and the first time that an analysis of the impact of such integration in terms of operational and financial aspects has been conducted. As the proposed model is a large-scale, mixed-integer linear problem, we developed a heuristic algorithm that, although making use of classical methods, works well for the intentions of this paper.

The remainder of this article is structured as follows. In Section \ref{sec:models}, the novel mathematical model that integrates the processes at the core of papermaking is proposed. In Section \ref{sec:solution}, the heuristic solution we developed is explained. In Section \ref{sec:experiments}, the tests conducted to analyse the behaviour of the proposed model are described. In Section \ref{sec:conclusions}, the main findings are detailed.

\section{Mathematical model of the three phases of the paper production process} \label{sec:models}

In this section, we present the mathematical model of the three phases of the paper production process in an integrated way. To amortise the uncertainties of demand, the rolling-horizon planning technique \citep{Sethi1991} was used. The long-term planning Phases 1 and 2 are solved for every period, and the short-term planning Phase 3 is solved for the sub-periods that compose the first period of the planning horizon.

Phase 1 is modelled based on a lot-sizing problem. The goal is to find the quantity of jumbos to be produced with a minimum total cost---constituted by production, machine setup and storage costs---in each period of the planning horizon and respecting the machine capacity of the factory.

To model Phase 2, we used a one-dimensional cutting-stock problem. The aim is to find the best way of cutting jumbos into reels for each period, minimising the costs of the material waste in the cutting process, the machine setup and the reel storage costs. It must meet the reel demand without exceeding the machine capacity.

Phase 3 is modelled based on a two-dimensional cutting-stock problem. We seek to find the best way to cut reels into sheets for each sub-period of the first period of the planning horizon, meeting the sheet demand and respecting the machine capacity. The costs considered in this phase are the material waste, the machine setup and the sheet storage costs.

The indexes, parameters and variables are defined next. We sought to maintain the same letter to specify notations referring to similar definitions throughout the model and indicated the phases they refer to using the number that follows. For example, $m_1$, $m_2$ and $m_3$ refer to the type of machine in phases 1, 2 and 3, respectively.\\

\noindent \textit{Indexes:} 

\noindent $t$: period in the planning horizon; varies in $t = 1, \ldots, T$, where $T$ is the total number of periods;

\noindent $\tau$: sub-period in the first period; varies in $\tau = 1, \ldots, \Theta$, where $\Theta$ is the total number of sub-periods in the first period;

\noindent $k$: paper grammage; varies in $k = 1, \ldots, K$, where $K$ is the quantity of types of grammage;

\noindent $m_1$: jumbo production machine type; varies in $m_1 = 1, \ldots, M_1$, where $M_1$ is the total number of different types of machines; machine $m_1$ makes jumbos with the length $L_{m_1}$;    

\noindent $m_2$: rewinder type; varies in $m_2 = 1, \ldots, M_2$, where $M_2$ is the total number of different types of rewinders;

\noindent $i_2$: reel type, with the dimensions $l_{i_2}$ for length and $w_{i_2}$ for width; varies in $i_2 = 1, \ldots, N\hspace{-0.12cm}f_2$, where $N\hspace{-0.12cm}f_2$ is the total number of different types of reels;    

\noindent $j_2$: one-dimensional cutting pattern; varies in $j_2 = 1, \ldots, N_{m_1}$, where $N_{m_1}$ is the pattern quantity for jumbos of type $m_1$;

\noindent $\{1,\ldots, {N\hspace{-0.12cm}f}_2\} =  S_2(1) \cup S_2(2) \cup \ldots\cup S_2(K)$, where $S_2(k)$=\{$i_2$ such that the grammage of reel $i_2$ is $k$\}.    

\noindent $m_3$: cutter type; varies in $m_3 = 1, \ldots, M_3$, where $M_3$ is the total number of different types of cutters;

\noindent $i_3$: sheet type, with the dimensions $l_{i_3}$ for length and $w_{i_3}$ for width; varies in $i_3 = 1, \ldots, {N\hspace{-0.12cm}f}_3$, where ${N\hspace{-0.12cm}f}_3$ is the total number of different types of sheets;

\noindent $j_3$: two-dimensional cutting pattern; varies in $j_3 = 1, \ldots, N_{i_2}$, where $N_{i_2}$ is the pattern quantity for reels of type $i_2$;

\noindent $\{1,\ldots, {N\hspace{-0.12cm}f}_3\} =  S_3(1) \cup S_3(2) \cup \ldots\cup S_3(K)$, where $S_3(k)$=\{$i_3$ such that the grammage of sheet $i_3$ is $k$\}; \\

\noindent \textit{Parameters:} 

\noindent $c_{1_{k, m_1, t}}$: production cost of jumbo of grammage $k$, in machine $m_1$, in period $t$, per weight unit;

\noindent $s_{1_{k, m_1, t}}$: setup cost for jumbo of grammage $k$, produced by machine $m_1$, in period $t$;

\noindent $h_{1_{k, t}}$: cost/kg for stocking jumbos of grammage $k$ in period $t$;

\noindent $b_{1_{k, m_1}}$: weight of jumbos of grammage $k$, produced by machine $m_1$;

\noindent $d_{1_{k, m_1, t}}$: demand of jumbos of grammage $k$ and length $L_{m_1}$ in period $t$; 

\noindent $f_{1_{k, m_1}}$: time for producing jumbos of grammage $k$ in machine $m_1$;

\noindent $g_{1_{k, m_1}}$: setup time for producing jumbos of grammage $k$ in machine $m_1$;

\noindent $C_{1_{t}}$: available time for jumbos production in period $t$;

\noindent $a^{j_2}_{2_{i_2, m_1}}$: quantity of reels of type $i_2$ cut from jumbo produced by machine $m_1$, according to pattern $j_2$. 

\noindent $c_{2_{k, t}}$: waste cost for paper of grammage $k$ during the jumbos cutting process, in period $t$, per length units;

\noindent $p^{j_2}_{2_{m_1}}$: paper waste in jumbos of type $m_1$ cutting, according to pattern $j_2$; calculated by $p^{j_2}_{2_{m_1}} = L_{m_1} - \sum_{k = 1}^{K}\sum_{i_2 \in S_2(k)} a^{j_2}_{2_{i_2, m_1}} l_{i_2}$;

\noindent $s^{j_2}_{2_{k, m_1, m_2, t}}$: setup cost for jumbo of grammage $k$, produced by machine $m_1$, cut by rewinder $m_2$, in period $t$, according to pattern $j_2$;

\noindent $h_{2_{i_2, t}}$: cost/kg for reels of type $i_2$ stock, in period $t$;

\noindent $b_{2_{i_2}}$: reels of type $i_2$ weight;

\noindent $d_{2_{i_2, t}}$: reels of type $i_2$ demand in period $t$; 

\noindent $\alpha_{i_2}$: increase of extra reels of type $i_2$ needed to supply Phase 3, for the periods $t > 1$; 

\noindent $\delta_{2_{i_2, t}}$: reels of type $i_2$ extra demand that must be supplied in period $t$ to feed Phase 3; as Phase 3 demand for $t > 1$ is not known, $\delta_{2_{i_2, t}} = \delta_{2_{i_2, 1}}$ for $t = 2, \ldots, T$ was considered;

\noindent $f^{j_2}_{2_{k, m_1, m_2}}$: time for cutting jumbos of grammage $k$, produced by machine $m_1$, cut in rewinder $m_2$, according to pattern $j_2$;

\noindent $g^{j_2}_{2_{k, m_1, m_2}}$: setup time for cutting jumbos of grammage $k$, produced by machine $m_1$, cut in rewinder $m_2$, according to pattern $j_2$;

\noindent $C_{2_t}$: available time for jumbos cutting in period $t$;

\noindent $a^{j_3}_{3_{i_3, i_2}}$: quantity of sheets of type $i_3$ cut from reels of type $i_2$, according to pattern $j_3$; 

\noindent $c_{3_{k, \tau}}$: waste cost for paper of grammage $k$ during reels cutting process, in sub-period $\tau$, per area units;

\noindent $p^{j_3}_{3_{i_2}}$: paper waste in reels of type $i_2$ cutting, according to pattern $j_3$; calculated by $p^{j_3}_{3_{i_2}} = l_{i_2} w_{i_2} - \sum_{k = 1}^{K}\sum_{i_3 \in S_3(k)} a^{j_3}_{3_{i_3, i_2}} l_{i_3} w_{i_3}$;     

\noindent $s^{j_3}_{3_{i_2, m_3, \tau}}$: setup cost for reels of type $i_2$ cutting, in cutter $m_3$, in sub-period $\tau$, according to pattern $j_3$; 

\noindent $h_{3_{i_3, \tau}}$: cost/kg for sheet of type $i_3$ stock, in sub-period $\tau$;

\noindent $b_{3_{i_3}}$: sheets of type $i_3$ weight;

\noindent $d_{3_{i_3, \tau}}$: sheets of type $i_3$ demand in sub-period $\tau$; 

\noindent $f^{j_3}_{3_{i_2, m_3}}$: time for reels of type $i_2$ cutting, in cutter $m_3$, according to pattern $j_3$;

\noindent $g^{j_3}_{3_{i_2, m_3}}$: setup time for reels of type $i_2$ cutting, in cutter $m_3$, according to pattern $j_3$;

\noindent $C_{3_{\tau}}$: available time for reels cutting in sub-period $\tau$;

\noindent $Q$: a sufficiently large number;\\

\noindent \textit{Variables:} 

\noindent $x_{1_{k, m_1, t}}$: quantity of jumbos of grammage $k$ produced by machine $m_1$, in period $t$;

\noindent $z_{1_{k, m_1, t}}$: binary variables that indicate whether jumbos of grammage $k$ were produced by machine $m_1$, in period $t$;   

\noindent $e_{1_{k, m_1, t}}$: quantity of jumbos of grammage $k$, produced by machine $m_1$ stocked in period $t$.

\noindent $y^{j_2}_{2_{k, m_1, m_2, t}}$: quantity of jumbos of grammage $k$, produced by machine $m_1$, cut in rewinder $m_2$, in period $t$, according to pattern $j_2$;

\noindent $z^{j_2}_{2_{k, m_1, m_2, t}}$: binary variables that indicate whether jumbos of grammage $k$, produced by machine $m_1$, were cut in rewinder $m_2$, in period $t$, according to pattern $j_2$;

\noindent $e_{2_{i_2, t}}$: reels of type $i_2$ stocked in period $t$. 

\noindent $y^{j_3}_{3_{i_2, m_3, \tau}}$: reels of type $i_2$ quantity cut in cutter $m_3$, in sub-period $\tau$, according to pattern $j_3$;

\noindent $z^{j_3}_{3_{i_2, m_3, \tau}}$: binary variables that indicate whether reels of type $i_2$ were cut in cutter $m_3$, in sub-period $\tau$, according to pattern $j_3$;

\noindent $e_{3_{i_3, \tau}}$: sheets of type $i_3$ quantity stocked in sub-period $\tau$.\\ 

The integrated paper production process is represented in Model \ref{eq:modelo}. As we integrate two already integrated models (Phases 1 and 2 together with Phases 2 and 3), this new formulation is termed as the Bi-Integrated Model.

\begin{subequations} \label{eq:modelo}
	\allowdisplaybreaks[1]
	\begin{align} 
	& \text{Minimise} \nonumber \\
	& \sum_{t=1}^{T} \sum_{m_1=1}^{M_1} \sum_{k=1}^{K} (c_{1_{k, m_1, t}} x_{1_{k, m_1, t}} + s_{1_{k, m_1, t}} z_{1_{k, m_1, t}} + h_{1_{k, t}} b_{1_{k, m_1}} e_{1_{k, m_1, t}}) + \label{eq:fase1_fo}\\    
	& \sum_{t=1}^{T} \sum_{m_1=1}^{M_1} \sum_{m_2=1}^{M_2} \sum_{k=1}^{K} \sum_{j_2=1}^{N_{m_1}} (c_{2_{k, t}} p^{j_2}_{2_{m_1}} y^{j_2}_{2_{k, m_1, m_2, t}} + s^{j_2}_{2_{k, m_1, m_2, t}} z^{j_2}_{2_{k, m_1, m_2, t}}) + \nonumber \\
	& \hspace{8cm} \sum_{t=1}^{T} \sum_{k=1}^{K} \sum_{i_2 \in S_2(k)} h_{2_{i_2, t}} b_{2_{i_2}} e_{2_{i_2, t}} + \label{eq:fase2_fo}\\    
	& \sum_{\tau=1}^{\Theta} \sum_{m_3=1}^{M_3} \sum_{k=1}^{K} \sum_{j_3=1}^{N_{i_2}} \sum_{i_2 \in S_2(k)} (c_{3_{k, \tau}} p^{j_3}_{3_{i_2}} y^{j_3}_{3_{i_2, m_3, \tau}} + s^{j_3}_{3_{i_2, m_3, \tau}} z^{j_3}_{3_{i_2, m_3, \tau}}) + \nonumber \\
	&\hspace{8cm}\sum_{\tau=1}^{\Theta} \sum_{k=1}^{K} \sum_{i_3 \in S_3(k)} h_{3_{i_3, \tau}} b_{3_{i_3}} e_{3_{i_3, \tau}} \label{eq:fase3_fo} \\
	& \text{subject to} \nonumber \\    
	& x_{1_{k, m_1, t}} + e_{1_{k, m_1, t - 1}} - e_{1_{k, m_1, t}} = d_{1_{k, m_1, t}} +  \sum_{m_2 = 1}^{M_2} \sum_{j_2 = 1}^{N_{m_1}} y^{j_2}_{2_{k, m_1, m_2, t}}, \quad k = 1, \ldots, K, \nonumber \\ 
	& \hspace{8cm}m_1 = 1, \ldots, M_1, \quad t = 1, \ldots, T.     \label{eq:integracao12} \\
	& \sum_{m_1=1}^{M_1} \sum_{k=1}^{K} (f_{1_{k, m_1}} x_{1_{k, m_1, t}} + g_{1_{k, m_1}} z_{1_{k, m_1, t}})  \leq C_{1_{t}}, \quad t = 1, \ldots, T; \label{eq:1_capacidade} \\
	& x_{1_{k, m_1, t}} \leq Q z_{1_{k, m_1, t}}, \quad k = 1, \ldots, K, \quad m_1 = 1, \ldots, M_1, \quad t = 1, \ldots, T; && \label{eq:1_preparacao}  \\
	&\sum_{m_1 = 1}^{M_1} \sum_{m_2 = 1} ^{M_2} \sum_{j_2 = 1}^{Nm_1} a^{j1}_{2_{i_2, m_1}} y^{j_2}_{2_{k, m_1, m_2, 1}} - e_{2_{i_2, 1}} = d_{2_{i_2, 1}} + \sum_{\tau = 1}^{\Theta} \sum_{m_3 = 1}^{M_3} \sum_{j_3 = 1}^{N_{i_2}} y^{j_3}_{3_{i_2, m_3, \tau}}, \nonumber \\
	& \hspace{8cm} i_2 \in S_2(k), \quad k = 1, \ldots, K. \label{eq:integracao23} \\
	& \sum_{m_1 = 1}^{M_1} \sum_{m_2 = 1}^{M_2} \sum_{j_2 = 1}^{N_{m_1}} \sum_{i_2 \in S_2(k)} (a^{j_2}_{2_{i_2, m_1}} y^{j_2}_{2_{k, m_1, m_2, t}} + e_{2_{i_2, t-1}} - e_{2_{i_2, t}}) = (1 + \alpha_{i_2}) d_{2_{i_2, t}}, \nonumber \\
	& \hspace{6cm} k = 1, \ldots, K, \quad i_2 \in S_2(k), \quad t = 2, \ldots, T; \label{eq:2_demanda_123} \\	
	& \sum_{m_1 = 1}^{M_1} \sum_{m_2 = 1}^{M_2} \sum_{k = 1}^{K} \sum_{j_2 = 1}^{N_{m_1}} (f^{j_2}_{2_{k, m_1, m_2}} y^{j_2}_{2_{k, m_1, m_2, t}} + g^{j_2}_{2_{k, m_1, m_2}} z^{j_2}_{2_{k, m_1, m_2, t}}) \leq C_{2_t}, \quad t = 1, \ldots, T; \label{eq:fase2_capacidade} \\
	&y^{j_2}_{2_{k, m_1, m_2, t}} \leq Q z^{j_2}_{2_{k, m_1, m_2, t}}, \quad k = 1, \ldots, K, \quad m_1 = 1, \ldots, M_1, \quad m_2 = 1, \ldots, M_2, \nonumber \\
	& \hspace{8cm} t = 1, \ldots, T, \quad j_2 = 1, \ldots, N_{m_1}; \label{eq:fase2_preparacao}  \\
	&\sum_{m_3 = 1}^{M_3} \sum_{j_3 = 1}^{N_{i_2}} \sum_{i_2 \in S(k)} (a^{j_3}_{3_{i_3, i_2}} y^{j_3}_{3_{i_2, m_3, \tau}} + e_{3_{i_3, \tau - 1}} - e_{3_{i_3, \tau}}) = d_{3_{i_3, \tau}}, \nonumber \\
	& \hspace{6cm} i_3 \in S_3(k), \quad k = 1, \ldots, K, \quad \tau = 1, \ldots, \Theta; \label{eq:fase3_demanda}\\
	&\sum_{m_3 = 1}^{M_3} \sum_{k = 1}^{K} \sum_{j_3 = 1}^{N_{i_2}} \sum_{i_2 \in S(k)} (f^{j_3}_{3_{i_2, m_3}} y^{j_3}_{3_{i_2, m_3, \tau}} + g^{j_3}_{3_{i_2, m_3}} z^{j_3}_{3_{i_2, m_3, \tau}}) \leq C_{3_{\tau}}, \quad \tau = 1, \ldots, \Theta; \label{eq:fase3_capacidade}\\
	&y^{j_3}_{3_{i_2, m_3, \tau}} \leq Q z^{j_3}_{3_{i_2, m_3, \tau}}, \quad j_3 = 1, \ldots, N_{i_2}, \quad i_2 \in S_2(k), \quad  k = 1, \ldots, K, \nonumber \\
	& \hspace{7cm} m_3 = 1, \ldots, M_3, \quad \tau = 1, \ldots, \Theta; \label{eq:fase3_preparacao} \\
	& x_{1_{k, m_1, t}}, \quad e_{1_{k, m_1, t}} \ge 0 \text{ and integers}, \quad z_{1_{k, m_1, t}} \in \{0, 1\}, \quad e_{1_{k, m_1, 0}} = 0, \nonumber \\
	& \hspace{5cm} k = 1, \ldots, K, \quad m_1 = 1, \ldots, M_1, \quad t = 1, \ldots, T. \label{eq:1_variaveis}    \\    
	& y^{j_2}_{2_{k, m_1, m_2, t}}, \quad e_{2_{i_2, t}} \ge 0 \text{ and integers}, \quad z^{j_2}_{2_{k, m_1, m_2, t}} \in \{0, 1\}, \quad e_{2_{i_2, 0}} = 0, \nonumber \\
	&k = 1, \ldots, K, \quad i_2 \in S_2(k), \quad  m_1 = 1, \ldots, M_1, \quad m_2 = 1, \ldots, M_2, \quad t = 1, \ldots, T, \nonumber \\
	& \hspace{10cm} j_2 = 1, \ldots, N_{m_1}.  \label{eq:fase2_variaveis} \\        
	& y^{j_3}_{3_{i_2, m_3, \tau}}, \quad e_{3_{i_3, \tau}} \ge 0 \text{ and integers}, \quad z^{j_3}_{3_{i_2, m_3, \tau}} \in \{0, 1\}, \quad e_{3_{i_3, 0}} = 0, \nonumber \\
	& k = 1, \ldots, K, \quad i_3 \in S_3(k), \quad i_2 \in S_2(k), \quad m_3 = 1, \ldots, M_3, \quad \tau = 1, \ldots, \Theta, \nonumber \\
	& \hspace{10cm} j_3 = 1, \ldots, N_{i_2}. \label{eq:fase3_variaveis}    
	\end{align}
\end{subequations} 

The objective function (\ref{eq:fase1_fo} + \ref{eq:fase2_fo} + \ref{eq:fase3_fo}) is a composition of the objective functions of the three phases. The portion (\ref{eq:fase1_fo}) minimises the production cost of jumbos, the machine setup and the stock costs (Phase 1). The portion (\ref{eq:fase2_fo}) refers to the costs of paper waste in the process of cutting jumbos into reels, rewinder setup and reel stock (Phase 2). The portion (\ref{eq:fase3_fo}) corresponds to the costs of paper waste when cutting reels into sheets, cutter setup and sheet stock (Phase 3).

The block (\ref{eq:integracao12}) to (\ref{eq:1_preparacao}) refers to Phase 1 constraints. Restriction (\ref{eq:integracao12}) is for meeting the total demand for jumbos, composed of the orders portfolio and the extra quantity needed to supply Phase 2. It is responsible for integrating Phases 1 and 2, since it contains variables from Phase 1---$x_{1_{k, m_1, t}}$ and $e_{1_{k, m_1, t}}$---and from Phase 2---$y^{j_2}_{2_{k, m_1, m_2, t}}$. (\ref{eq:1_capacidade}) guarantees that the machine capacity is not exceeded for each period. In (\ref{eq:1_preparacao}), setup binary variables are set to 1 when there is jumbo production and to 0 otherwise.

The block (\ref{eq:integracao23}) to (\ref{eq:fase2_preparacao}) refers to Phase 2 constraints. Restrictions (\ref{eq:integracao23}) and (\ref{eq:2_demanda_123}) are for meeting the demand for reels. It is necessary to furnish the requested reels as final products and the quantity needed to supply Phase 3. For the first period, this is guaranteed by Restriction \ref{eq:integracao23}, that integrates Phases 2---through the variables $y^{j_2}_{2_{k, m_1, m_2, 1}}$ and $e_{2_{i_2, 1}}$--- and Phase 3---through the variables $y^{j_3}_{3_{i_2, m_3, \tau}}$. As there is information about the demand of sheets only for the sub-periods of the first period of the planning horizon, for the remainder periods ($t>1$) we considered an increase $\alpha_{i_2}$ in the quantity of reels to be supplied to Phase 3. (\ref{eq:fase2_capacidade}) guarantees that the capacity of the rewinders is not exceeded for each period. In (\ref{eq:fase2_preparacao}), setup binary variables are set to 1 if jumbos are cut and to 0 otherwise.

The block (\ref{eq:fase3_demanda}) to (\ref{eq:fase3_preparacao}) refers to Phase 3 constraints. Restriction (\ref{eq:fase3_demanda}) is for meeting the demand for sheets. (\ref{eq:fase3_capacidade}) guarantees that the capacity of the cutters is not exceeded for each sub-period. (\ref{eq:fase3_preparacao}) makes setup binary variables set to 1 whenever reels are cut and to 0 otherwise. 

The block (\ref{eq:1_variaveis}) to (\ref{eq:fase3_variaveis}) refers to all the constraints on the variables. (\ref{eq:1_variaveis}) to (\ref{eq:fase3_variaveis}) restrict the variables $x_{1_{k, m_1, t}}$, $e_{1_{k, m_1, t}}$, $y^{j_2}_{2_{k, m_1, m_2, t}}$, $e_{2_{i_2, t}}$, $y^{j_3}_{3_{i_2, m_3, \tau}}$ and $e_{3_{i_3, \tau}}$ to non-negative integer values and variables $z_{1_{k, m_1, t}}$, $z^{j_2}_{2_{k, m_1, m_2, t}}$ and $z^{j_3}_{3_{i_2, m_3, \tau}}$ to binary values. The initial jumbo, reel and sheet stocks are set to 0.

Model \ref{eq:modelo} represents a large-scale, mixed-integer linear programming problem. Therefore, it becomes challenging to find optimal solutions using conventional exact methods for instances of significant dimensions. Different techniques can be explored to solve this novel model. Following, we present the method we developed in detail.

\section{Solution method} \label{sec:solution}

The first difficulty in solving Model \ref{eq:modelo} is the existence of binary machine setup variables in the three phases of the production process. The strategy adopted was not explicitly to consider the costs and setup times in the model, but rather to add them empirically to costs and production times, avoiding that the solution obtained by Model \ref{eq:modelo} becomes infeasible in practice. Therefore, some parameters were updated to take into account this increase due to machine setup. Calling $c^{'}_{1_{k, m_1, t}}$ and $f^{'}_{1_{k, m_1}}$ the new cost and production time, respectively, for Phase 1; $c^{'}_{2_{k, t}}$ and $f^{'j_2}_{2_{k, m_1, m_2}}$ the new paper waste cost and cutting time, respectively, for Phase 2; and $c^{'}_{3_{k, \tau}}$ and $f^{'j_3}_{3_{i_2, m_3}}$ the new paper waste cost and cutting time, respectively, for Phase 3, the updated objective function is defined by Expression (\ref{eq:fo_new}).

\begin{multline}
\sum_{t=1}^{T} \sum_{m_1=1}^{M_1} \sum_{k=1}^{K} (c^{'}_{1_{k, m_1, t}} x_{1_{k, m_1, t}} + h_{1_{k, t}} b_{1_{k, m_1}} e_{1_{k, m_1, t}}) + \\
\sum_{t=1}^{T} \sum_{m_1=1}^{M_1} \sum_{m_2=1}^{M_2} \sum_{k=1}^{K} \sum_{j_2=1}^{N_{m_1}} c^{'}_{2_{k, t}} p^{j_2}_{2_{m_1}} y^{j_2}_{2_{k, m_1, m_2, t}} + 
\sum_{t=1}^{T} \sum_{k=1}^{K} \sum_{i_2 \in S_2(k)} h_{2_{i_2, t}} b_{2_{i_2}} e_{2_{i_2, t}} + \\
\sum_{\tau=1}^{\Theta} \sum_{m_3=1}^{M_3} \sum_{k=1}^{K} \sum_{j_3=1}^{N_{i_2}} \sum_{i_2 \in S_2(k)} c^{'}_{3_{k, \tau}} p^{j_3}_{3_{i_2}} y^{j_3}_{3_{i_2, m_3, \tau}} + \sum_{\tau=1}^{\Theta} \sum_{k=1}^{K} \sum_{i_3 \in S_3(k)} h_{3_{i_3, \tau}} b_{3_{i_3}} e_{3_{i_3, \tau}} \label{eq:fo_new}
\end{multline}

The capacity restrictions also had to be updated, according to Expressions (\ref{eq:1_capacidade_relaxado}), (\ref{eq:2_capacidade_relaxado}) and (\ref{eq:3_capacidade_relaxado}) for Phases 1, 2 and 3, respectively.

\begin{subequations} \label{eq:capacity_new}
	\allowdisplaybreaks[1]
	\begin{align}
	& \sum_{m_1=1}^{M_1} \sum_{k=1}^{K} f^{'}_{1_{k, m_1}} x_{1_{k, m_1, t}} \leq C_{1_{t}}, \quad t = 1, \ldots, T; \label{eq:1_capacidade_relaxado} \\
	& \sum_{m_1 = 1}^{M_1} \sum_{m_2 = 1}^{M_2} \sum_{k = 1}^{K} \sum_{j_2 = 1}^{N_{m_1}} f^{'j_2}_{2_{k, m_1, m_2}} y^{j_2}_{2_{k, m_1, m_2, t}} \leq C_{2_t}, \quad t = 1, \ldots, T; \label{eq:2_capacidade_relaxado} \\
	&\sum_{m_3 = 1}^{M_3} \sum_{k = 1}^{K} \sum_{j_3 = 1}^{N_{i_2}} \sum_{i_2 \in S(k)} f^{'j_3}_{3_{i_2, m_3}} y^{j_3}_{3_{i_2, m_3, \tau}} \leq C_{3_{\tau}}, \quad \tau = 1, \ldots, \Theta; \label{eq:3_capacidade_relaxado}
	\end{align}
\end{subequations}

The second difficulty, that all the variables must be integers, is bypassed by solving the relaxed mathematical model and, with the approximated solution, applying a rounding heuristic to come up with an integer-valued solution. The relaxed mathematical model is then composed of Expression (\ref{eq:fo_new}), as the objective function, and the restrictions (\ref{eq:integracao12}), (\ref{eq:1_capacidade_relaxado}), (\ref{eq:integracao23}), (\ref{eq:2_capacidade_relaxado}), (\ref{eq:fase3_demanda}), (\ref{eq:3_capacidade_relaxado}), (\ref{eq:1_variaveis}) (with $x_{1_{k, m_1, t}}$ and $e_{1_{k, m_1, t}}$ relaxed), (\ref{eq:fase2_variaveis}) (with $y^{j_2}_{2_{k, m_1, m_2, t}}$ and $e_{2_{i_2, t}}$ relaxed) and (\ref{eq:fase3_variaveis}) (with $y^{j_3}_{3_{i_2, m_3, \tau}}$ and $e_{3_{i_3, \tau}}$ relaxed).

The last difficult, the large number of possible cutting patterns, is overcome by applying the Simplex Method with column generation \citep{Gomory1961, Gomory1965} to the relaxed mathematical model. 

\begin{figure}[!ht]
	\centering
	\includegraphics[width=\textwidth]{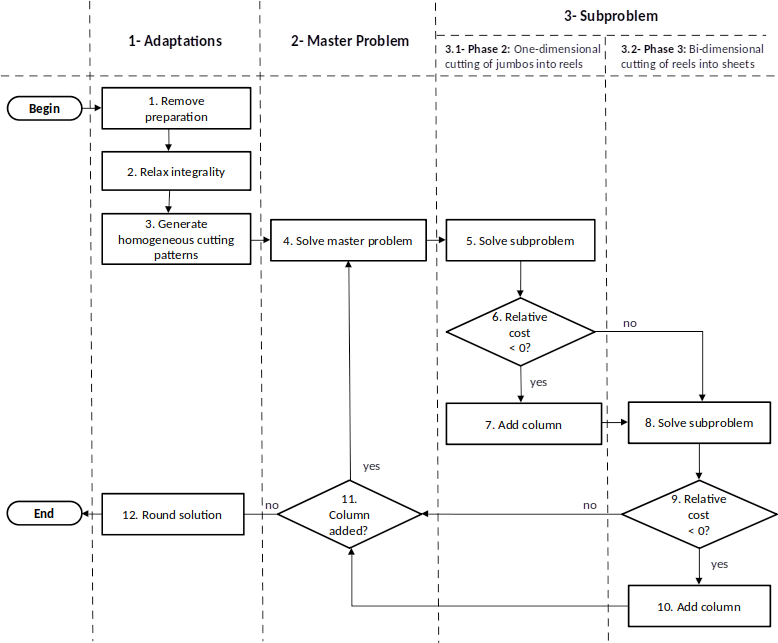}
	\caption{Solution heuristic method flowchart.}
	\label{figura:metodo}
\end{figure}

The flowchart in Figure \ref{figura:metodo} illustrates the proposed solution method. In the adaptations stage, the setup variables are removed and the integer variables are relaxed. In addition, homogeneous cutting patterns are generated for Phases 2 and 3, resulting in the initial restricted master problem. 

In the master problem stage, the master problem is solved and the optimal values of the dual variables $\pi$ are obtained.

The last stage, cutting-patterns generation, is divided into two parts, one to generate one-dimensional patterns for the jumbos and another to generate two-dimensional patterns for the reels. In the first case, the sub-problems represented by Model \ref{passo5} are solved, one for each combination of grammage $k$, machine type $m_1$ and period $t$, to obtain the best patterns $a^{j_2}_{2_{i_2, m_1}}$:

\begin{subequations} \label{passo5}
	\begin{align}
	&\text{Maximise} \nonumber \\
	&\sum_{i_2 \in S_2(k)} (c^{'}_{2_{k, t}} l_{2_{i_2}} + \pi_{(K M_1 T + T + (t - 1) {N\hspace{-0.05cm}f}_2 + \sum_{i = 1}^{k - 1} (S_2(i) + i_2))}) a^{j_2}_{2_{i_2, m_1}} \\
	&\text{subject to}\nonumber \\
	& \sum_{i_2 \in S_2(k)} l_{2_{i_2}} a^{j_2}_{2_{i_2, m_1}} \leq L_{m_1}, && \\
	& a^{j_2}_{2_{i_2, m_1}} \geq 0 \text{ and integer}, \quad k = 1, \ldots, K, \quad m_1 = 1, \ldots, M_1, \quad t = 1, \ldots, T, \quad   && 
	\end{align}
\end{subequations}        

If the relative cost $\hat{c}^{j_2}_{2_{k, m_1, t}}$ (Equation \ref{passo6}) of the new column is negative, it is inserted into the master problem. This is applied to every new generated column.

\begin{equation} \label{passo6}
\begin{split}
\hat{c}^{j_2}_{2_{k, m_1, t}} = & c^{'}_{2_{k, t}} L_{m_1} - \sum_{i_2 \in S_2(k)} c^{'}_{2_{k, t}} l_{2_{i_2}} a^{j_2}_{2_{i_2, m_1}} - \pi_{((k - 1) M_1 T + (m_1 - 1) T + t)} - \\
&- \pi_{(K M_1 T + T + (t - 1) {N\hspace{-0.05cm}f}_2 + \sum_{i = 1}^{k - 1} (S_2(i) + i_2))} a^{j_2}_{2_{i_2, m_1}} - \\
&- \sum_{m_2 = 1}^{M_2} \pi_{(K M_1 T + T + T {N\hspace{-0.05cm}f}_2 + t)} f^{'j_2}_{2_{k, m_1, m_2}}
\end{split}
\end{equation}        

In any event, the algorithm follows to the second part of this stage and solves several other sub-problems, as presented in Section \ref{sec:phase3_subproblems}, to obtain the best two-dimensional cutting patterns $\textbf{a}^{j_3}_{3_{i_2}}$, one for each sub-period combination $\tau$ and reel type $i_2$. Again, if the relative cost $\hat{c}^{j_3}_{3_{i_2, \tau}}$ (Equation \ref{passo9}) of the new column for each generated pattern is negative, it is added to the master problem. In any event, the algorithm returns to the master problem stage.

\begin{equation} \label{passo9}
\begin{split}
\hat{c}^{j_3}_{3_{i_2, \tau}} = & c^{'}_{3_{k, \tau}} l_{2_{i_2}} w_{2_{i_2}} - \sum_{i_3 \in S_3(k)} c^{'}_{3_{k, t}} l_{3_{i_3}} w_{3_{i_3}} a^{j_3}_{3_{i_3, i_2}} - \pi_{(K M_1 T + T + i_2)} \\            
&- \sum_{i_3 \in S_3(k)} \pi_{(K M_1 T + T + T {N\hspace{-0.12cm}f}_2 + T + (\tau - 1) {N\hspace{-0.12cm}f}_3 + \sum_{i = 1}^{k - 1} (S_3(i) + i_3))} a^{j_3}_{3_{i_3, i_2}} - \\
&- \sum_{m_3 = 1}^{M_3} \pi_{(K M_1 T + T + T {N\hspace{-0.12cm}f}_2 + T + \Theta {N\hspace{-0.12cm}f}_3 + \tau)} f^{'j_3}_{3_{i_2, m_3}}  \\
\end{split}
\end{equation}        

If any new column was added during the steps of these sub-problems, the master problem is solved again and a new iteration of the cutting-patterns generation stage is started. Otherwise, the current relaxed solution can no longer be improved and the adaptations step recurs to round the solution to integer values using a heuristic procedure. 

\subsection{Sub-problems of Phase 3} \label{sec:phase3_subproblems}

The sub-problems of Phase 3 are solved according to \citet{Gomory1965}, with two-stage guillotine cutting patterns. First, for each reel $i_2$ of dimensions $w_{i_2} \times l_{i_2}$, we define the best cutting patterns for every possible strip of fixed width $w_2 = w_{i_2}$ and varied lengths $l_{i_3}$, $i_3 \in S_3(k)$. Second, we determine the number of times each strip must be used to maximise the utility function of the reel.

Let $L_{i_3^*}$ be the set of paper sheets that can be cut from the strip $w_2 \times l_{i_3^*}$, $i_3^* \in S_3(k)$; that means $L_{i_3^*} = \{i_3 | l_{i_3} \le l_{i_3^*}, i_3 \in S_3(k)\}$ if paper trimming is allowed and $L_{i_3^*} = \{i_3 | l_{i_3} = l_{i_3^*}, i_3 \in S_3(k)\}$ otherwise. To find the best cutting pattern $\boldsymbol{\alpha^{j_3}_{3_{i_2, i_3^*}}}$ for the strip $w_2 \times l_{i_3^*}$, Model \ref{eq:strip} must be solved. $\alpha^{j_3}_{3_{i_2, i_3^*, i_3}}$ is the number of paper sheets of type $i_3$ in the strip $w_2 \times l_{i_3^*}$ and $V_{i_2, i_3^*}$ is the utility value of the strip.

\begin{subequations}  \label{eq:strip}
	\begin{align}
	& V_{i_2, i_3^*} = \text{maximise} && \nonumber \\
	&\sum_{i_3 \in L_{i_3^*}} (c^{'}_{3_{k, t}} l_{3_{i_3}} w_{i_3} + \pi_{(K M_1 T + T + T {N\hspace{-0.12cm}f}_2 + T + (\tau - 1) {N\hspace{-0.12cm}f}_3 + \sum_{i = 1}^{k - 1} (S_3(i) + i_3))}) \alpha^{j_3}_{3_{i_2, i_3^*, i_3}} \\
	& \text{subject to} && \nonumber \\    
	& \sum_{i_3 \in L_{i_3^*}} w_{3_{i_3}} \alpha^{j_3}_{3_{i_2, i_3^*, i_3}} \leq w_2, && \\
	& \alpha^{j_3}_{3_{i_2, i_3^*, i_3}} \geq 0 \text{ and integer}, \quad i_3 \in L_{i_3^*}. && 
	\end{align}
\end{subequations}    

Let $r_k$ be the number of different lengths $l_{i_3}$, $i_3 \in S_3(k)$. After calculating the best patterns for all the $r_k$ strips of reel $i_2$, we solve Model \ref{mod:passo2_subproblema_fase3} to find $\beta^{j_3}_{3_{i_2, i_3^*}}$, the number of times each strip $w_2 \times l_{i_3^*}$, $i_3^* \in S_3(k)$ must be used.

\begin{subequations} \label{mod:passo2_subproblema_fase3}
	\begin{align}
	& \text{Maximise} && \nonumber \\
	& \sum_{i_3^* = 1}^{r_k} V_{i_2, i_3^*} \beta^{j_3}_{3_{i_2, i_3^*}} && \\ 
	& \text{subject to} && \nonumber \\    
	& \sum_{i_3^* = 1}^{r_k} \beta^{j_3}_{3_{i_2, i_3^*}} \le l_{i_2}, && \\
	& \beta^{j_3}_{3_{i_2, i_3^*}} \ge 0 \text{ and integer}, \quad i_2 \in S_2(k), \quad k = 1, \ldots, K.  && 
	\end{align}
\end{subequations}

The final cutting pattern $\boldsymbol{a^{j_3}_{3_{i2}}}$ is obtained by Equation \ref{eq:a3}.

\begin{equation} \label{eq:a3}
\boldsymbol{a^{j_3}_{3_{i2}}} = \left( \begin{array}{c}
a^{j_3}_{3_{i2, 1}} \\
\ldots \\
a^{j_3}_{3_{i2, {N\hspace{-0.12cm}f}_3}} \end{array} \right) = \left( \begin{array}{c}
\sum_{i_3^* = 1}^{r_k} \alpha^{j_3}_{3_{i_2, i_3^*, 1}} \beta^{j_3}_{3_{i_2, i_3^*}} \\
\ldots \\
\sum_{i_3^* = 1}^{r_k} \alpha^{j_3}_{3_{i_2, i_3^*, {N\hspace{-0.12cm}f}_3}} \beta^{j_3}_{3_{i_2, i_3^*}} \end{array} \right)
\end{equation}

\subsection{Solution rounding}

\begin{figure}
	\centering
	\includegraphics[width=\textwidth]{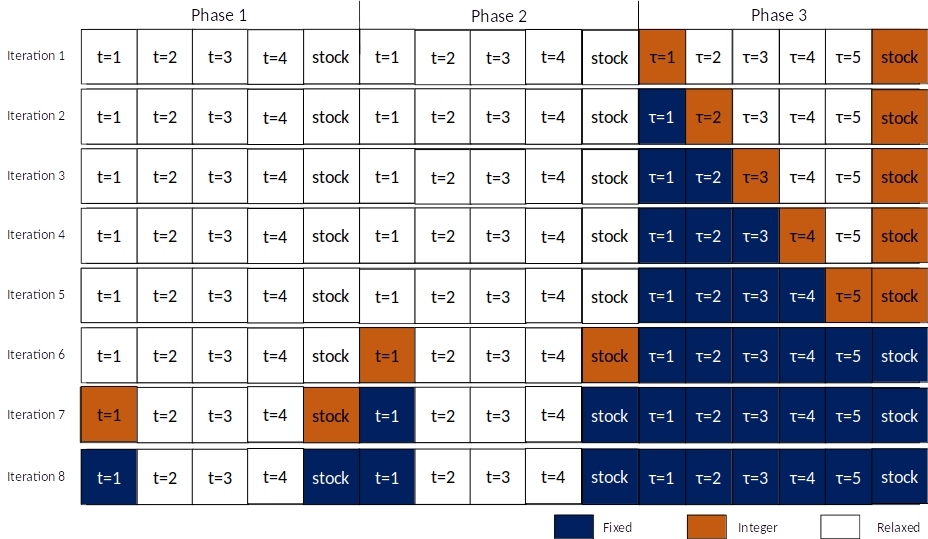}	
	\caption{Relax-and-fix.}
	\label{figura:esquema_arredondamento}
\end{figure}

As the model generated by the heuristic method is large-scale, with dimensions varying according to the number of different types of products in each phase of the process, it is not possible to solve it with integer variables in a reasonable time, as the preliminary experiments confirmed. Therefore, we adopted a technique to round the solution called \textit{relax-and-fix} \citep{Wolsey1998}. This consists of dividing the variables into sub-sets and progressively rearranging them into subproblems across the planning horizon, similar to what is successfully done by \citet{Brito2014}. Figure \ref{figura:esquema_arredondamento} illustrates its application to the Bi-Integrated Model rounding for an example with four periods, represented by $t$, and five sub-periods, $\tau$. It starts from the relaxed solution and executes a procedure of repetitively fixing the previous (sub-) period and obtaining the integer solution of the current (sub-) period, while keeping the remaining variables not yet rounded as real numbers. It proceeds to the variables of the third phase and then passes to those of the second phase and finishes with the first phase.

\section{Computational experiments} \label{sec:experiments}

We used Matlab and IBM CPLEX v. 12.6.1 (student license) to implement the solution method. For the rounding procedure, CPLEX solution timeout was configured to 20 minutes for each iteration. After that, if the solver was not able to find the optimum, it returns the best solution found so far. This timeout was rarely achieved, with the median total time to round the solution not passing the tens of seconds for the set of computational experiments executed.

Table \ref{tab:classes} presents the 24 classes used in the tests, with 20 instances each. The data were generated with random variables based on a local paper industry in Brazil and on \citet{Poltroniere2008a} to simulate the sensitive information (such as the costs) that was not provided by the enterprise.

\begin{table}[!h]
	\scriptsize
	\centering
	\tbl{Configurations of classes}
	{\begin{tabular}{cccccccc} \toprule
			\textbf{\#} & \begin{tabular}[c]{@{}l@{}}${N\hspace{-0.05cm}f}_2$\\${N\hspace{-0.05cm}f}_3$\end{tabular} & \textbf{\begin{tabular}[c]{@{}l@{}}Stock\\ cost\end{tabular}} & \textbf{\begin{tabular}[c]{@{}l@{}}Paper\\ trimming\\ allowed \end{tabular}} & $M_1$ & $M_2$ & $M_3$ & \textbf{\begin{tabular}[c]{@{}l@{}}\# \\work\\ shifts\end{tabular}} \\ \midrule
			1  & 5         & Normal           & Yes   & 3  & 3  & 2  & 1              \\
			2  & 5         & Normal           & No   & 3  & 3  & 2  & 1              \\
			3  & 5         & Normal           & Yes   & 3  & 3  & 2  & 2              \\
			4  & 5         & Normal           & No   & 3  & 3  & 2  & 2              \\
			5  & 5         & Normal           & Yes   & 3  & 3  & 2  & 3              \\
			6  & 5         & Normal           & No   & 3  & 3  & 2  & 3              \\
			7  & 5         & High             & Yes   & 3  & 3  & 2  & 1              \\
			8  & 5         & High             & No   & 3  & 3  & 2  & 1              \\
			9  & 5         & High             & Yes   & 3  & 3  & 2  & 2              \\
			10 & 5         & High             & No   & 3  & 3  & 2  & 2              \\
			11 & 5         & High             & Yes   & 3  & 3  & 2  & 3              \\
			12 & 5         & High             & No   & 3  & 3  & 2  & 3              \\
			13 & 9         & Normal           & Yes   & 6  & 6  & 4  & 1              \\
			14 & 9         & Normal           & No   & 6  & 6  & 4  & 1              \\ 
			15 & 9         & Normal           & Yes   & 6  & 6  & 4  & 2              \\ 
			16 & 9         & Normal           & No   & 6  & 6  & 4  & 2              \\  
			17 & 9         & Normal           & Yes   & 6  & 6  & 4  & 3              \\ 
			18 & 9         & Normal           & No   & 6  & 6  & 4  & 3              \\   
			19 & 9         & High             & Yes   & 6  & 6  & 4  & 1              \\ 
			20 & 9         & High             & No   & 6  & 6  & 4  & 1              \\ 
			21 & 9         & High             & Yes   & 6  & 6  & 4  & 2              \\ 
			22 & 9         & High             & No   & 6  & 6  & 4  & 2              \\ 
			23 & 9         & High             & Yes   & 6  & 6  & 4  & 3              \\ 
			24 & 9         & High             & No   & 6  & 6  & 4  & 3              \\  \bottomrule
	\end{tabular}}
	\label{tab:classes}
\end{table}

The number of work shifts is used to calculate the capacity of each type of machine. The change in the number of turns is used to vary the productive capacity of the plant.

The parameter intervals/definitions are summarised next. It is important to note that the number of types of a given parameter is different from its allowed range. For example, the number of types of paper grammage is set to 1 for all classes; however, its value can be varied for each instance in the interval [35; 300] $g/m^2$. Data for all the cases are available at \href{https://github.com/amandaortega/bi-integrated-model}{github.com/amandaortega/bi-integrated-model}. \\

\noindent \textit{General Parameters:} 

\noindent The diameter of jumbos and reels varies in the interval [300; 500] cm;

\noindent The paper thickness varies in the interval [190; 250] u.m.;

\noindent The grammage varies in the interval [35; 300] $g/m^2$;

\noindent $K$ (grammage number) $= 1$. \\

\noindent \textit{Phase 1:}

\noindent $c^{'}_{1_{k, m_1, t}}$ varies in the interval R\$ [0.015; 0.025]/kg;

\noindent $h_{1_{k, t}}$: normal cost varies in the interval R\$ [0.0075; 0.0125]/kg; high cost varies in the interval R\$ [0.009 ; 0.015 ]/kg;

\noindent $L_{m_1}$ varies in the interval [10; 20] m; 

\noindent $b_{1_{k, m_1}} = L_{m_1} \times \frac{\pi}{thickness} \times \frac{diameter^2}{4} \times grammage$;

\noindent $d_{1_{k, m_1, t}} = 0$ (no jumbos demanded as final product); 

\noindent $f^{'}_{1_{k, m_1}}$  varies in the interval [30; 60] min;

\noindent $C_{1_{t}} = $ 8h $\times$ \#work shifts $\times M_1 \times \theta$. \\

\noindent \textit{Phase 2:}

\noindent $c^{'}_{2_{k, t}}$ = R\$ $\frac{\sum_{m_1 = 1}^{M_1} c_{1_{k, m_1, t}}}{50000 M_1} $ / $cm^2$;

\noindent $h_{2_{i_2, t}} = 0,5 h_{1_{k, t}}$ / kg;

\noindent $l_2{_{i_2}}$  varies in the interval [3; 9] m;

\noindent $b_{2_{i_2}} = l_2{_{i_2}} \times \frac{\pi}{thickness} \times \frac{diameter^2}{4} \times grammage$;

\noindent $d_{2_{i_2, t}}$  varies in the interval [0; 100]; 

\noindent $f^{'j_2}_{2_{k, m_1, m_2}}$  varies in the interval [30; 60] min;

\noindent $C_{2_t}$  $ = $ 8h $\times$ \#work shifts $\times M_2 \times \theta$. \\

\noindent \textit{Phase 3:}

\noindent $c^{'}_{3_{k, \tau}} = $ R\$ $\frac{1.5 \sum_{m_1 = 1}^{M_1} \sum_{t = 1}^{T} c_{1_{k, m_1, t}}}{10000 M_1 T} $ / $cm^2$;

\noindent $h_{3_{i_3, \tau}} = 0.5 \frac{\sum_{k = 1}^{K} \sum_{t = 1}^{T} h_{1_{k, t}}}{K T}$ / kg; 

\noindent $l_{3_{i_3}}$ varies in the interval [30; 100] cm;

\noindent $w_{3_{i_3}}$ varies in the interval [30; 100] cm;

\noindent $b_{3_{i_3}} = l_3{_{i_2}} \times w_{3_{i_3}} \times grammage$;

\noindent $d_{3_{i_3, \tau}}$ varies in the interval [0; 500]; 

\noindent $f^{'j_3}_{3_{i_2, m_3}}$ varies in the interval [20; 40] min;

\noindent $C_{3_{\tau}} = $ 8h $\times$ \#work shifts $\times M_3$. \\

For the classes with more than one work shift, an increase of 20\% per additional work shift was considered in the jumbo production costs, $c^{'}_{1_{k, m_1, t}}$, and reel and sheet cutting costs, $c^{'}_{2_{k, t}}$ and $c^{'}_{3_{k, \tau}}$, respectively, due to the extra labor needed.

In the following sections, the results presented in figures use the median statistic measure because it is less sensitive to extreme values. When showing the results in tables, however, both mean and standard deviation measures are given.

\subsection{Comparison of the Bi-Integrated Model with other strategies to solve the problem}

\begin{figure}[!h]
	\centering
	\includegraphics[width=0.6\textwidth]{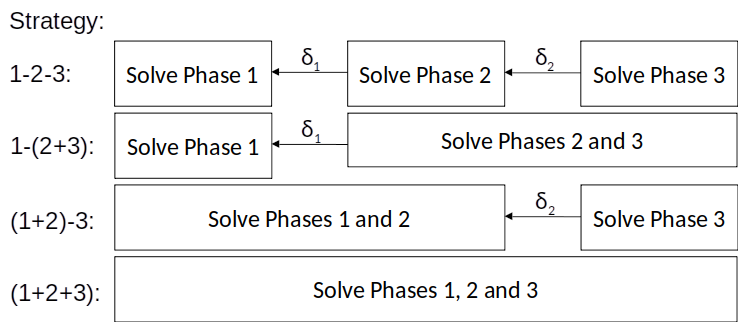}
	\caption{Strategies to solve the paper production problem.}
	\label{fig:strategies}
\end{figure}

Because the Bi-Integrated Model proposed in this paper is the pioneer dealing with the three phases of paper production, there is no work in the literature with which we can compare the results. Therefore, in this section, we compare the Bi-Integrated Model with other strategies to solve the problem, as illustrated in Figure \ref{fig:strategies}. 

Strategy 1-2-3 consists of considering each phase of the process separately. The solution flow follows the opposite direction. That is, Phase 3 is solved and then the resulting number of reels, $\delta_2$, is calculated. The demand for the reel portfolio is increased with this extra value and Phase 2 is solved. Based on the obtained result, the number of jumbos required, denoted by $\delta_1$, is calculated and added to the quantity in the portfolio of this type of product. Finally, Phase 1 is solved.

Strategy 1-(2+3) consists of solving Phases 2 and 3 in an integrated model and calculating the quantity $\delta_1$ of jumbos needed to supply these phases. With this extra value added to the demand, we finally solve Phase 1.

Strategy (1+2)-3 consists of integrating Phases 1 and 2 of the process. Phase 3 is solved and the quantity $\delta_2$ is calculated. This value is then added to the reel portfolio and Phases 1 and 2 are solved in an integrated model.

Strategy (1+2+3) is the strategy proposed by the Bi-Integrated Model.

In this section, the results obtained by class 1 for the four strategies are analysed according to several criteria. There is no particular reason why class 1 was chosen, and the analysis did not differ significantly for the remaining classes.

The parameter $\alpha_{i_2}$, the increase of the demand of reels of type $i_2$ to supply Phase 3 for the periods $t>1$ is calculated by Equation (\ref{eq:alpha}), for all the strategies.

\begin{equation}
\alpha_{i_2} = \frac{d_{2_{i_2, 1}} + \delta_2}{d_{2_{i_2, 1}}} - 1, \quad i_2 \in S_2(k), \quad k = 1, \ldots, K. \label{eq:alpha}
\end{equation}

\subsubsection{Costs}

\begin{figure}[!h]
	\centering
	\subfigure[Total costs.]{
		{\includegraphics[width=.45\textwidth]{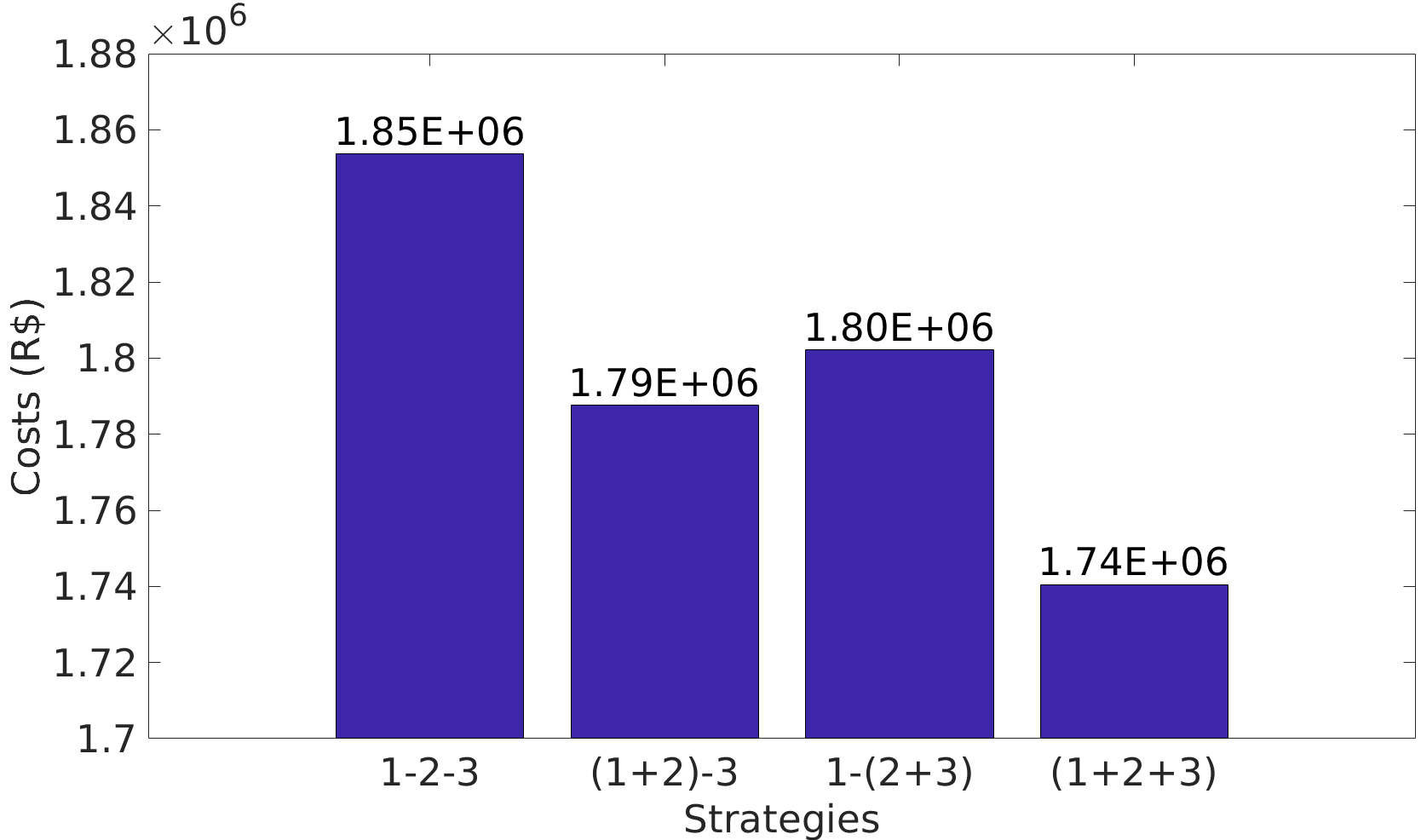}}\label{fig:custos_totais}}
	\subfigure[Gap between relaxed optimum  and rounded solutions.]{
		{\includegraphics[width=.45\textwidth]{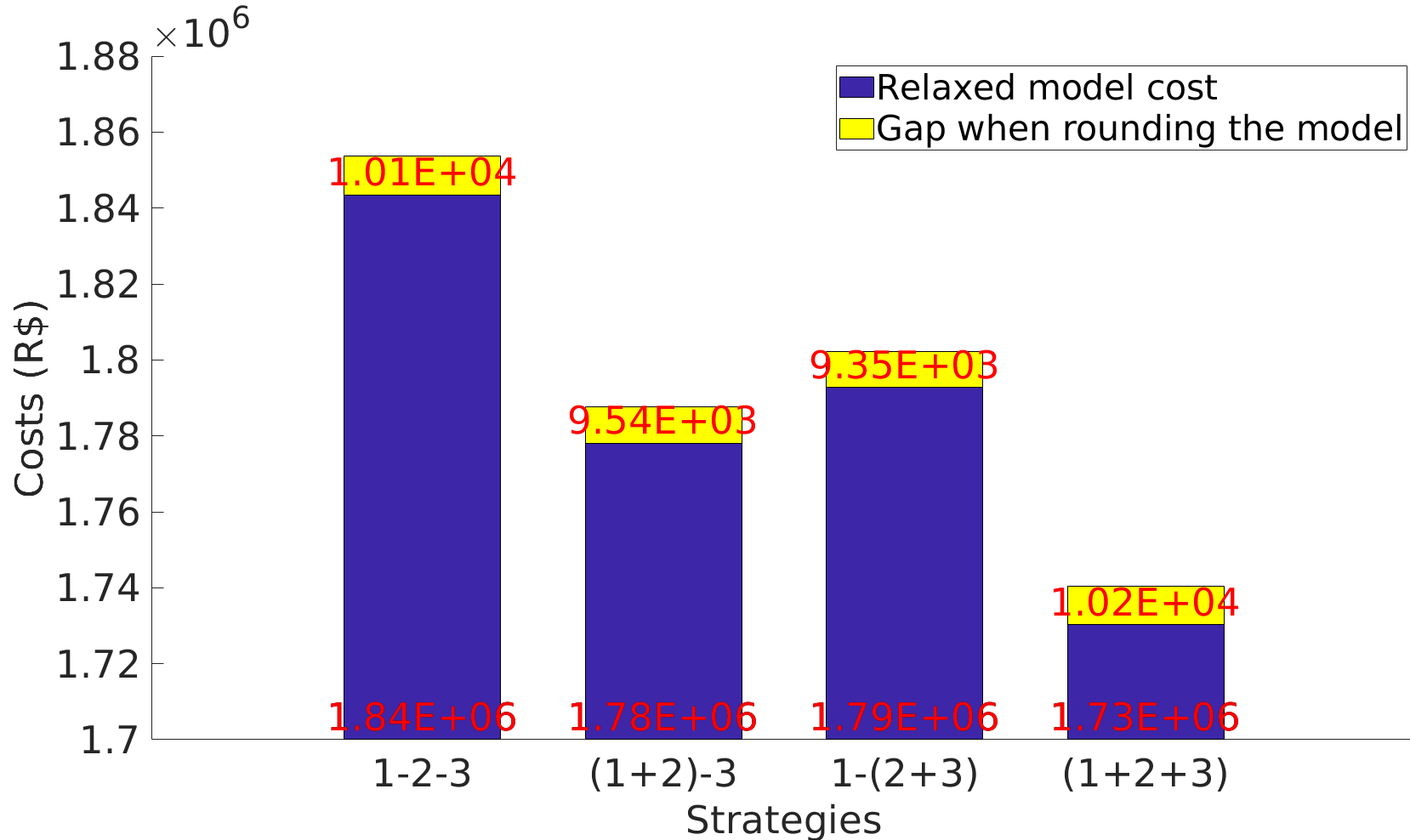}}\label{fig:diferenca_custos}}
	\caption{Total costs among the four strategies.} \label{fig:custos}
\end{figure}

Figure \ref{fig:custos} compares the total costs for the four strategies. Figure \ref{fig:custos_totais} shows a cost decrease as the process steps are integrated, reaching the value of nearly 6\% when moving from Strategy 1-2-3 to Strategy (1+2+3). In Figure \ref{fig:diferenca_custos}, the shares of the solutions due to the rounding step are highlighted. For all strategies, the maximum increment obtained was 0.59\% above the relaxed value, which indicates good performance of the rounding algorithm.

\begin{figure}[!h]
	\centering
	\subfigure[Total costs per phase.]{
		{\includegraphics[width=.6\textwidth]{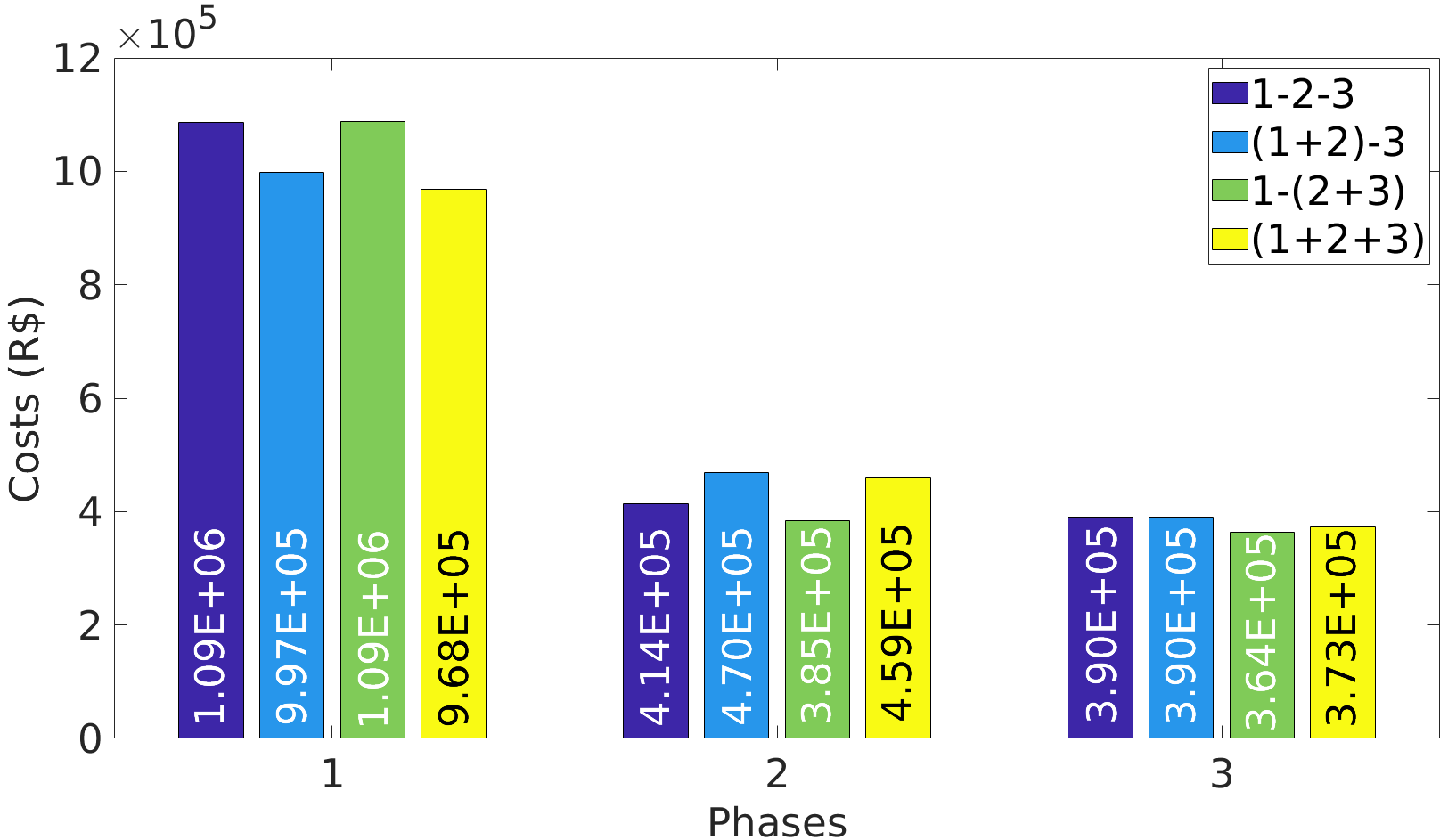}}\label{fig:custos_fase}}
	\subfigure[Production/cutting costs per phase.]{
		{\includegraphics[width=.45\textwidth]{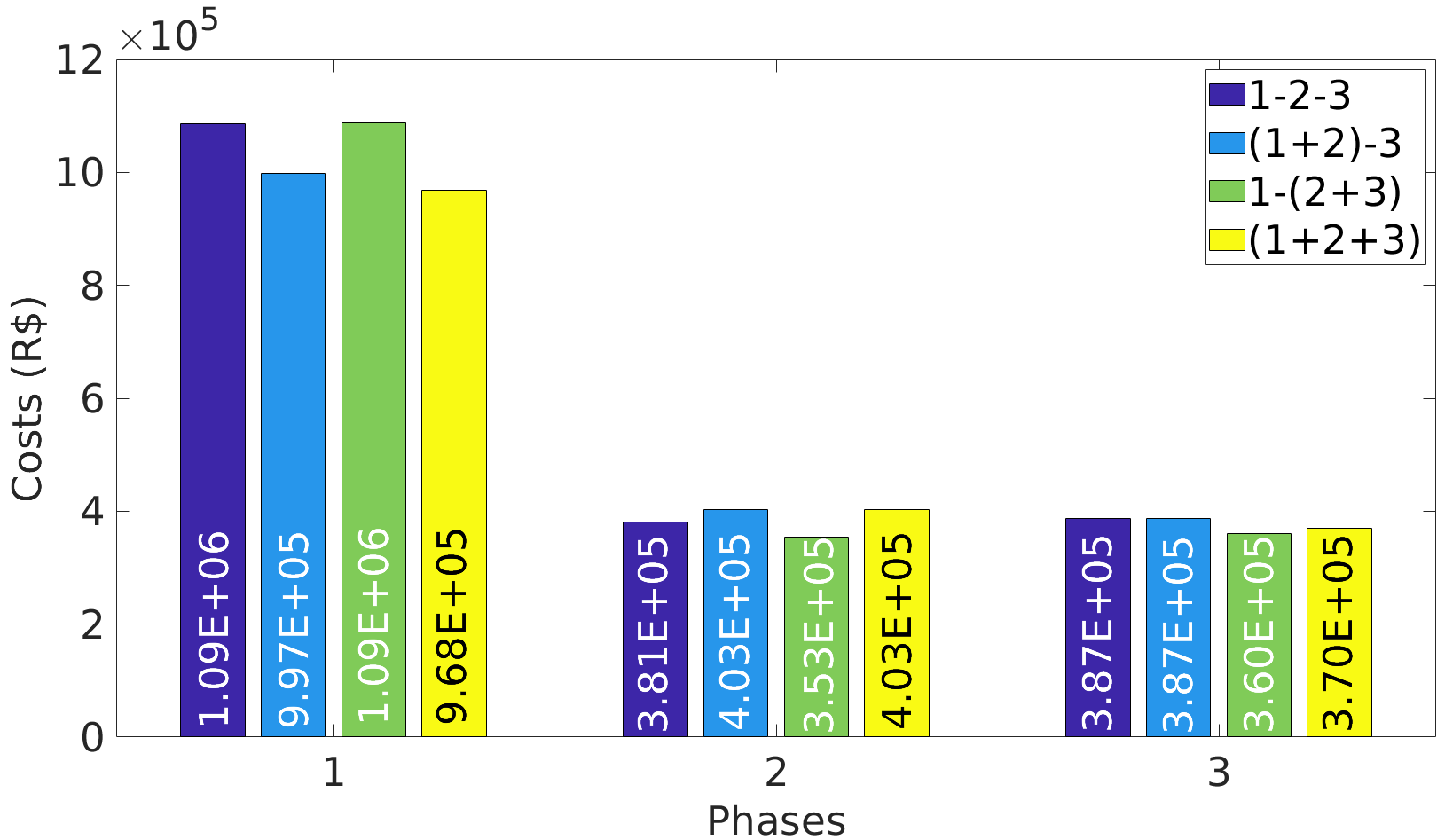}}\label{fig:custos_producao}}
	\subfigure[Stock costs per phase.]{
		{\includegraphics[width=.45\textwidth]{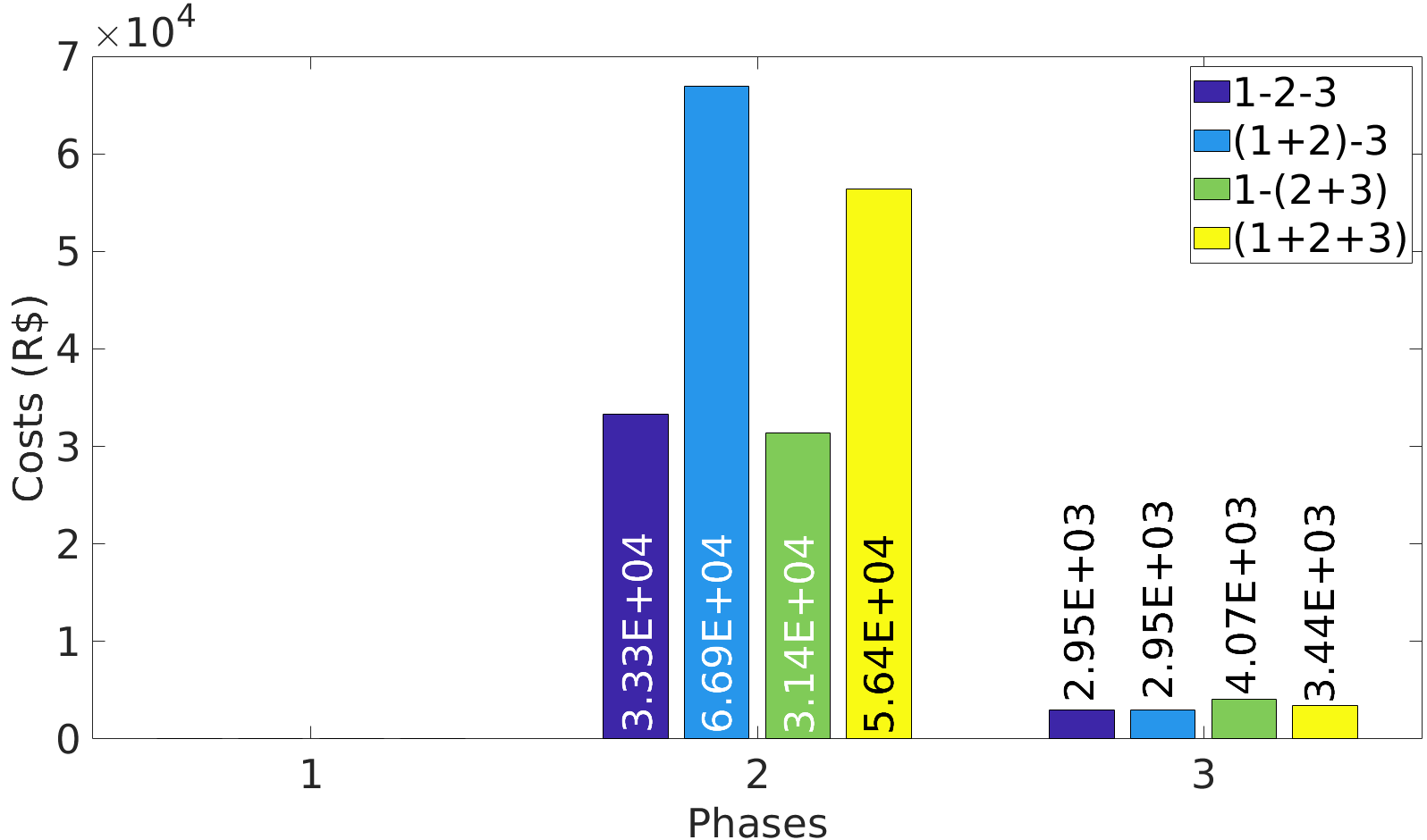}}\label{fig:custos_estoque}}    
	\caption{Total costs per phase among the four strategies.} \label{fig:custos_fase_todas}
\end{figure}

Figure \ref{fig:custos_fase_todas} shows the costs distribution per phase among the four strategies. By integrating phases, the result of one phase often worsens, while for other phases it improves. For example, in Figure \ref{fig:custos_fase}, by integrating Phases 1 and 2, Phase 2 costs are increased, but Phase 1 costs are decreased. As the drop in the cost of Phase 1 is more significant than the increase in the cost of Phase 2, the total cost decreased. Another example is what happens when integrating the three phases; when compared to Strategy 1-(2+3), Phase 2 and Phase 3 costs have increased but, as there is a greater reduction of Phase 1 costs, the total cost turns out to be lower. 

In Figures \ref{fig:custos_producao} and \ref{fig:custos_estoque}, the costs per phase separated by category are presented. Figure \ref{fig:custos_producao} shows production costs for Phase 1 and cutting costs for Phases 2 and 3; Figure \ref{fig:custos_estoque} shows stock costs for all phases. Primarily noticeable is that, regardless of the strategy adopted, there is no stock cost for Phase 1, since it is costly (financially and operationally) to stock jumbos. Another important point is that Phase 3 costs do not vary substantially from one strategy to another, either for cut or stock. Phase 2 costs show a slightly higher variation, caused mainly by the inventory. Nevertheless, what makes a difference in the final solution is the production of jumbos in Phase 1. This analysis should convince the reader that good planning of the process as a whole, with the entire chain taken into account, makes the total costs drop considerably. This is the concept of the Bi-Integrated Model proposed in this work.

\subsubsection{Material waste} \label{sec:analise_desperdicio}

\begin{figure}[!h]
	\centering
	\includegraphics[width=0.6\textwidth]{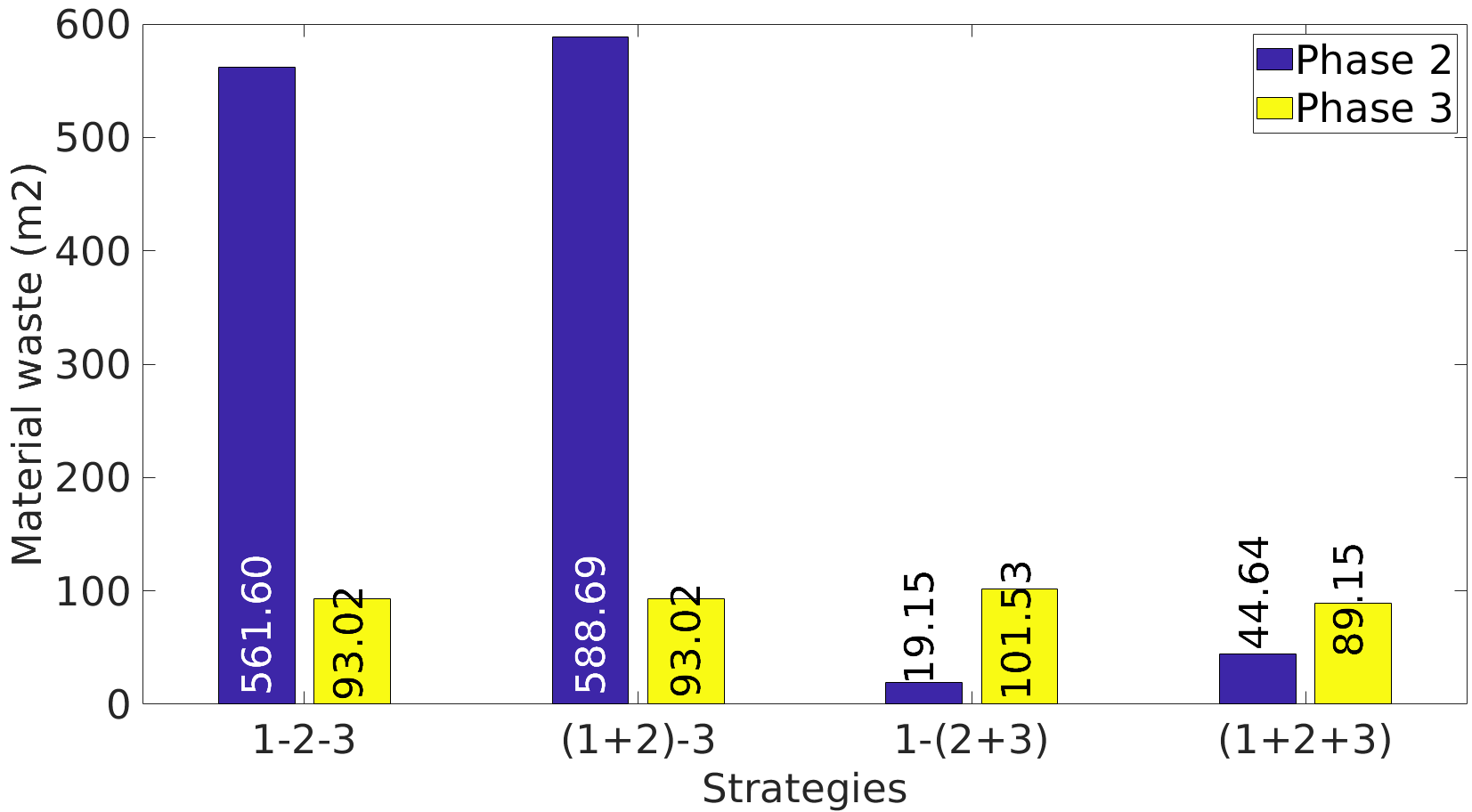}
	\caption{Material waste per phase among the four strategies.}
	\label{fig:desperdicio}
\end{figure}

Figure \ref{fig:desperdicio} shows the material waste for Phases 2 and 3 among the four strategies. Strategies 1-2-3 and (1+2)-3 obtain extremely high values of material waste for Phase 2, because the Phase 3 solution does not take into account the best combinations of reels that, cut from jumbos in Phase 2, will result in less waste. This yields cutting patterns with a small loss in Phase 3, but, when converted to extra reel demand in Phase 2, significant waste is generated. By integrating the two cutting stages, Phase 2 waste decreases considerably.

\subsubsection{Stock of products} \label{sec:analise_estoque}

\begin{figure}[!h]
	\centering
	\includegraphics[width=0.6\textwidth]{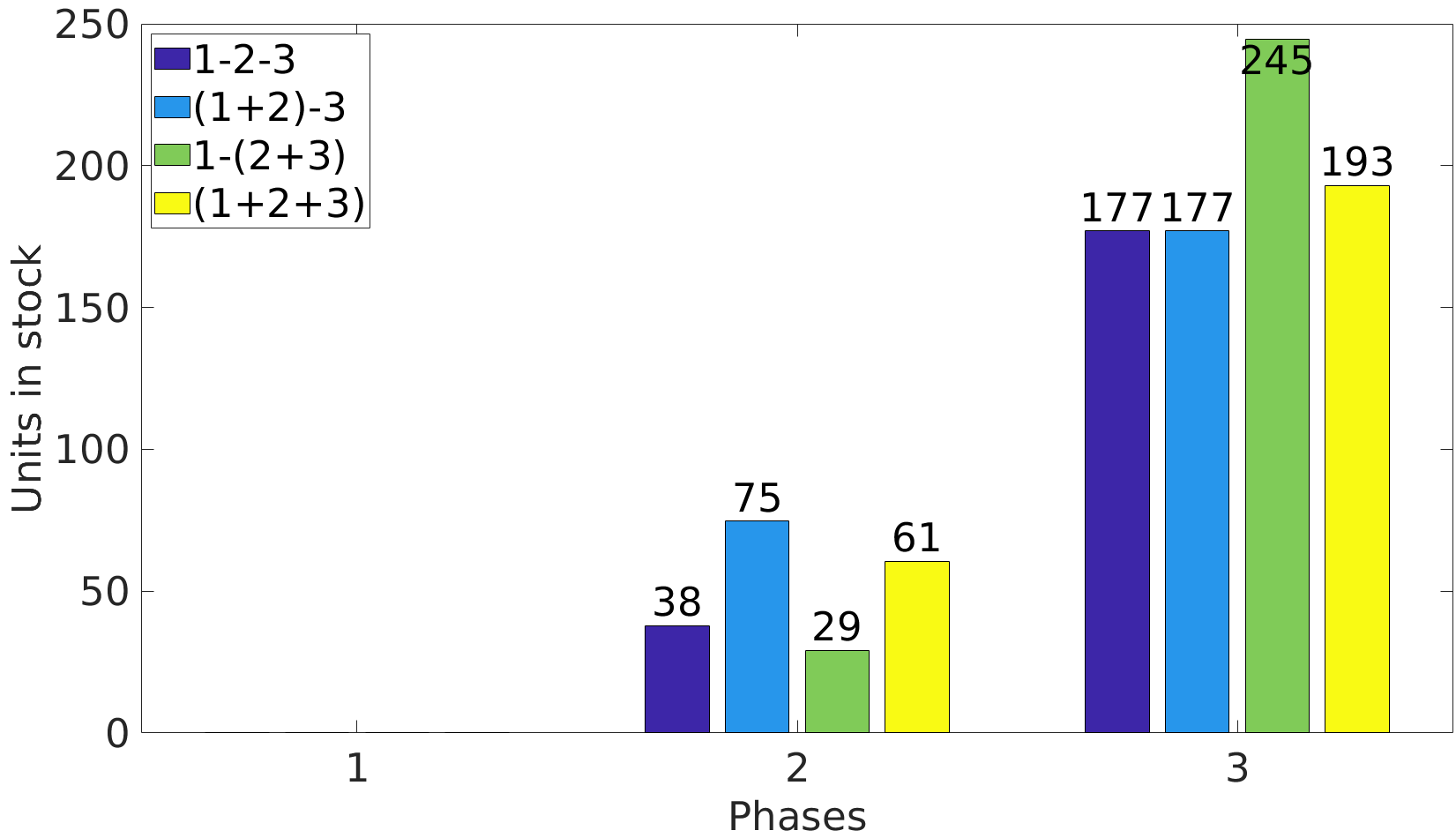}
	\caption{Stock per period/sub-period.}
	\label{fig:estoque}
\end{figure}

Figure \ref{fig:estoque} shows the number of units stocked per (sub-) period. As anticipated in Figure \ref{fig:custos_estoque}, there is no jumbo stock in any of the strategies. In addition, there is a tendency to use more stock when integrating the production phases, as better cutting patterns are used as a result of anticipating the cutting of some items in the planning horizon.

\subsubsection{Machine capacity} \label{sec:analise_capacidade}

\begin{figure}[!h]
	\centering
	\includegraphics[width=0.6\textwidth]{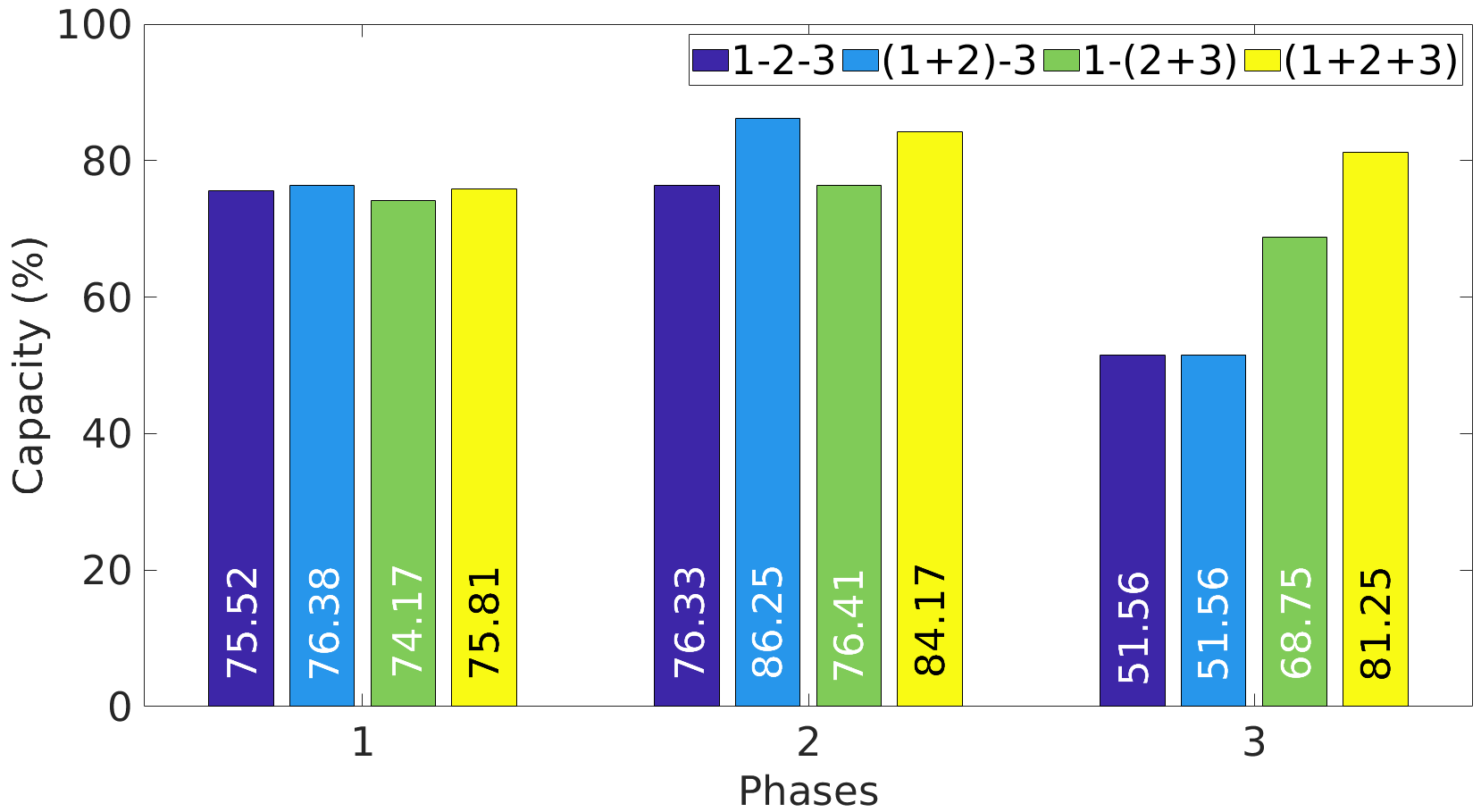}
	\caption{Maximum capacity used per phase among the four strategies.}
	\label{fig:capacidade}
\end{figure}

Figure \ref{fig:capacidade} shows the average values of maximum machine capacity per phase among the four strategies. Note there is a tendency to obtain higher maximum-capacity values as the phases are integrated. As better cutting patterns are obtained as a result of anticipating items from later periods, more machine-processing capacity is demanded.

\subsubsection{Number of columns} \label{sec:analise_colunas}

\begin{figure}[!h]
	\centering
	\subfigure[Composition of the models columns.]{
		{\includegraphics[width=.45\textwidth]{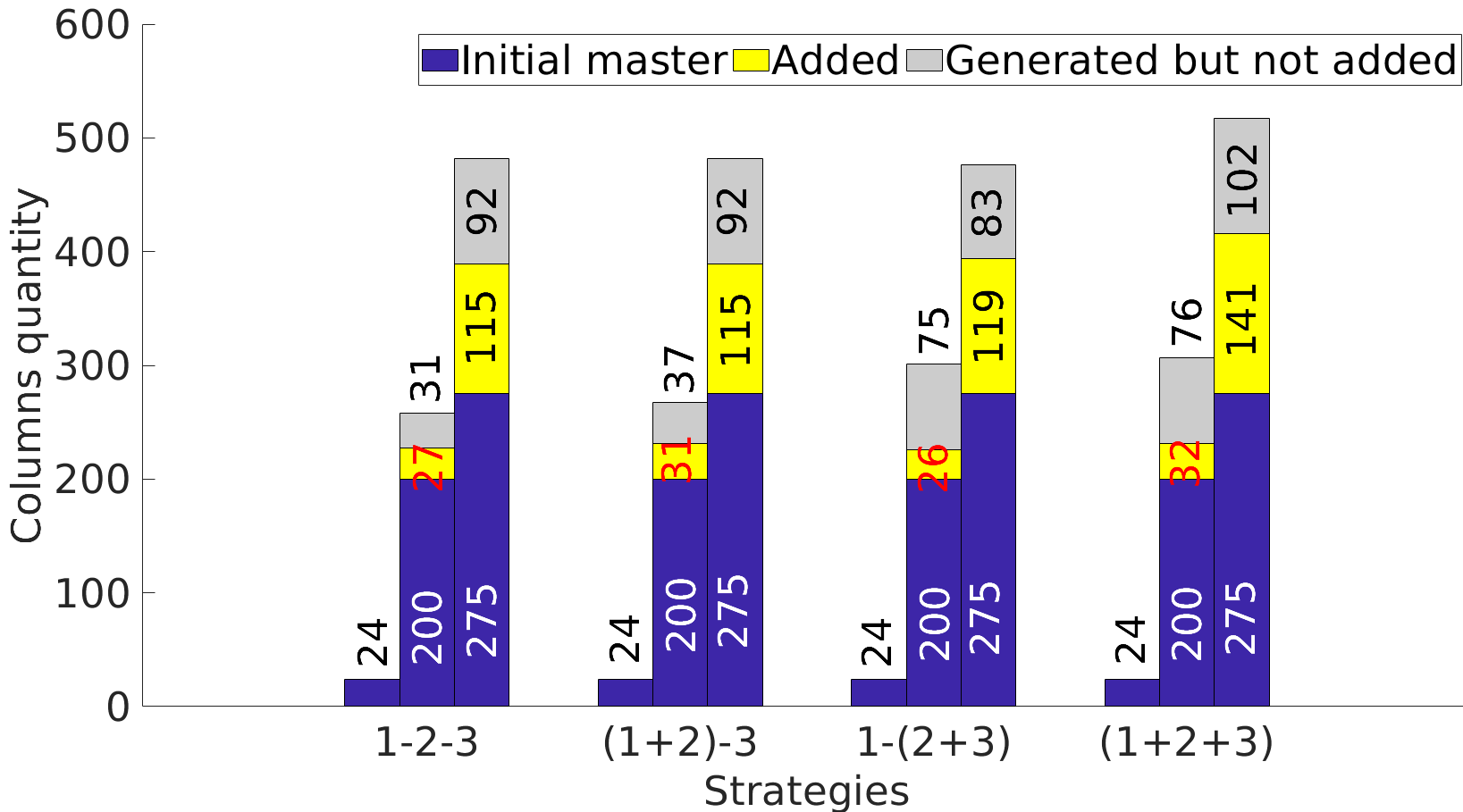}}\label{fig:composicao_colunas}}
	\subfigure[Composition of the columns used in the solution.]{
		{\includegraphics[width=.45\textwidth]{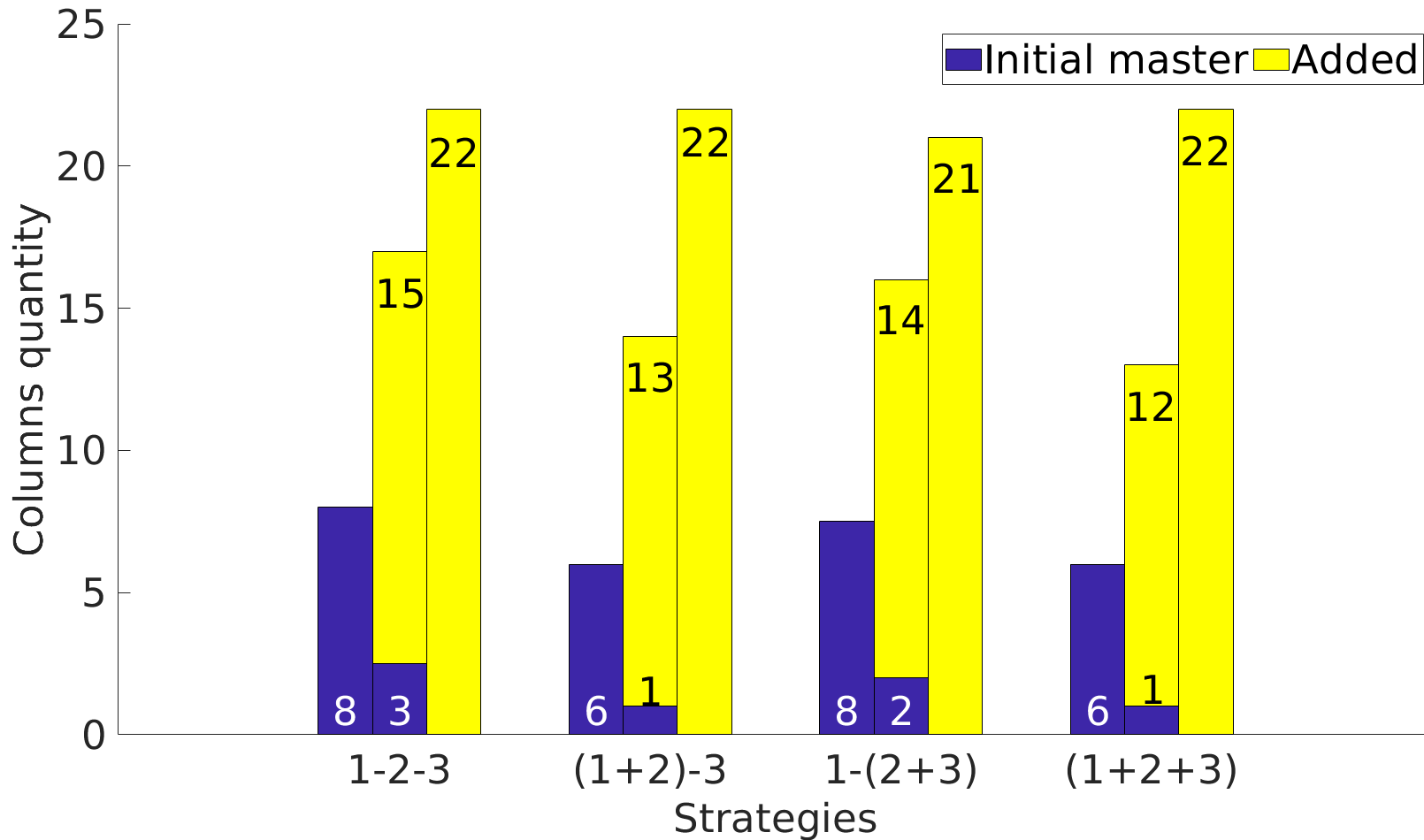}}\label{fig:colunas_usadas}}
	\caption{Analysis of the number of columns, by phase and among the four strategies.} \label{fig:analise_colunas}
\end{figure}

Figure \ref{fig:analise_colunas} shows an analysis of the quantities and origin of the columns used in the models for the three phases and among the four strategies. In Figure \ref{fig:composicao_colunas}, for Phases 2 and 3 there is a clear increasing trend of columns generated and inserted when integrating phases of the process, since new possibilities of patterns appear, in such a way as to obtain a lower global optimum. It is no coincidence that Strategy (1+2+3) is the one with the most significant number of columns generated and inserted, in the two cutting phases.

One may observe that, in the linear relaxation, the number of columns used in the optimal solution is equal to the size of the base matrix of each strategy (columns in blue in Figure \ref{fig:composicao_colunas}). In the rounded solution, however, the number of columns used in the model solution decreases, as shown in Figure \ref{fig:colunas_usadas}. This also indicates the origin of the columns, that is, whether they were part of the initial master problem or added by the column generation procedure. There is a tendency for a decrease in the columns of the initial master problem for Phases 2 and 3, as they are composed of homogeneous cutting patterns and present larger amounts of waste. Another important finding is that the number of columns used in the solution remains approximately constant in all strategies. What makes the solution better, therefore, is not the use of a greater number of columns, but rather, the substitution of some patterns for better ones. This fact is positive because the frequency with which the cutting knives should be repositioned due to the change of pattern remains constant.

\subsubsection{Processing time} \label{sec:analise_tempos}

\begin{figure}[!h]
	\centering
	\includegraphics[width=0.6\textwidth]{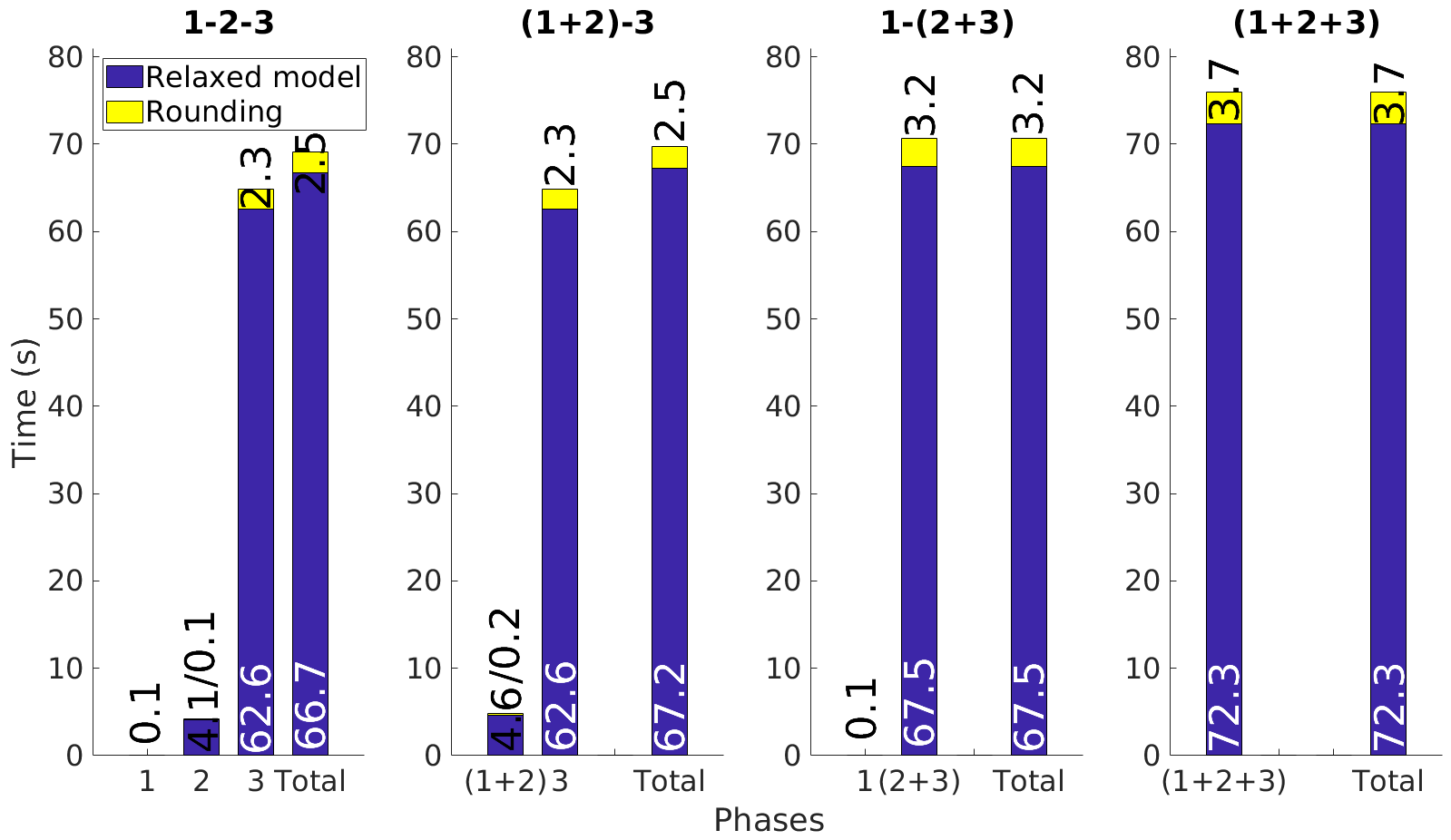}
	\caption{Heuristic processing time per phase among the four strategies.}
	\label{fig:tempos}
\end{figure}

Figure \ref{fig:tempos} shows the median processing time, in seconds, obtained by the proposed heuristic method by phase among the four strategies. Most of the processing time lies in solving the relaxed model, although the rounding procedure is not negligible. Although the complexity of the models grows considerably when integrating phases of the process, the same does not happen with the processing time, which only undergoes a small increase---units of seconds---, making the use of such strategies computationally possible . 

\subsection{Analysis of the Bi-Integrated Model in different parameter settings}

This sub-section presents the results of the 24 test classes, analysing several relevant factors. The results presented here are restricted to Strategy (1+2+3), but all the solutions were obtained for the four strategies and can be found online at \href{https://github.com/amandaortega/bi-integrated-model}{github.com/amandaortega/bi-integrated-model}. \\

\subsubsection{Total cost earnings} \label{sec:resumo_diminuicao_custos}

Table \ref{tab:diferenca_custos_estrategias} contains the results (mean and standard deviation) of the percentage difference in cost among the solutions obtained by the Bi-Integrated Model (Strategy (1+2+3)) compared to the other strategies, considering the relaxed and the rounded solutions. Analysing the standard deviation between them for classes 1--12 and 13--24, one may note more stable results achieved by the relaxed model. The rounding procedure brings high variability to the solutions. Another remarkable point is that, given the mean and the standard deviation for each relaxed value in Table \ref{tab:diferenca_custos_estrategias}, one can infer that Strategy (1+2+3) outperforms the other strategies in 100\% of the instances for all classes. 

In terms of the rounded solution, Strategy (1+2+3) obtained considerable reductions in the total cost, reaching the maximum value of 26.65\% of decline compared to Strategy (1-2-3) for class 18 and outstanding performances in general. However, based on the obtained values of the standard deviation, one can infer the existence of some instances where the results of Strategy (1+2+3) did not outperform the others for some particular cases. In the solution rounding stage, it may occur, for the other strategies, that the obtained integer quantities in each phase result in a smaller overall cost, due to the nature of the rounding process. 

\begin{table}[!h]
	\scriptsize
	\centering
	\caption{Percentage difference in total cost between Strategy (1+2+3) and the others}
	\label{tab:diferenca_custos_estrategias}
	\begin{tabular}{lrrrrrr} \hline
		\multirow{2}{*}{\#} & \multicolumn{2}{c}{$\Delta \frac{(1+2+3)}{1-2-3}$ (\%)} & \multicolumn{2}{c}{$\Delta \frac{(1+2+3)}{(1+2)-3}$ (\%)} & \multicolumn{2}{c}{$\Delta \frac{(1+2+3)}{1-(2+3)}$ (\%)} \\ 
		& Relaxed & Rounded & Relaxed & Rounded & Relaxed & Rounded \\
		\hline
		1 & -6.8$\pm$4.8 & -6.1$\pm$11.4 & -2.4$\pm$2.3 & -1.9$\pm$10.9 & -5.5$\pm$5.2 & -6.1$\pm$7.4 \\
		2 & -6.5$\pm$4.0 & -14.4$\pm$11.4 & -2.5$\pm$2.0 & -11.4$\pm$11.3 & -4.8$\pm$4.3 & -7.5$\pm$7.6 \\
		3 & -7.0$\pm$4.7 & -9.7$\pm$8.9 & -2.4$\pm$2.3 & -5.8$\pm$7.8 & -6.3$\pm$5.7 & -13.1$\pm$12.6 \\
		4 & -6.5$\pm$4.2 & -16.3$\pm$15.8 & -2.4$\pm$1.9 & -13.5$\pm$16.1 & -5.3$\pm$4.6 & -6.7$\pm$14.0 \\
		5 & -7.0$\pm$4.7 & -1.3$\pm$46.8 & -2.4$\pm$2.2 & -0.1$\pm$53.3 & -6.2$\pm$5.7 & 3.1$\pm$46.3 \\
		6 & -6.5$\pm$4.4 & -18.2$\pm$23.2 & -2.4$\pm$2.0 & -15.2$\pm$23.3 & -5.3$\pm$4.4 & -9.3$\pm$11.9 \\
		7 & -6.6$\pm$4.4 & -7.6$\pm$4.5 & -2.5$\pm$2.5 & -3.5$\pm$3.6 & -5.2$\pm$4.8 & -7.9$\pm$7.1 \\
		8 & -6.1$\pm$3.8 & -12.8$\pm$9.8 & -2.5$\pm$2.1 & -9.9$\pm$9.8 & -4.4$\pm$3.9 & -7.6$\pm$5.2 \\
		9 & -6.8$\pm$4.7 & -9.3$\pm$9.4 & -2.5$\pm$2.4 & -9.3$\pm$20.5 & -5.5$\pm$5.1 & -10.3$\pm$12.1 \\
		10 & -6.5$\pm$4.1 & -21.6$\pm$14.3 & -2.5$\pm$2.0 & -19.0$\pm$14.9 & -4.9$\pm$4.4 & -8.6$\pm$7.9 \\
		11 & -6.9$\pm$4.6 & -10.0$\pm$13.2 & -2.4$\pm$2.4 & -5.7$\pm$14.4 & -6.3$\pm$5.7 & -9.0$\pm$13.0 \\
		12 & -6.5$\pm$4.2 & -12.7$\pm$34.7 & -2.4$\pm$2.0 & -9.8$\pm$35.6 & -5.3$\pm$4.6 & -4.3$\pm$25.4 \\ \hline
		Mean & -6.6$\pm$0.3 & -11.7$\pm$12.6 & -2.4$\pm$0.2 & -8.8$\pm$14.1 & -5.4$\pm$0.6 & -7.3$\pm$11.7 \\ \hline
		13 & -10.2$\pm$3.5 & -12.3$\pm$5.2 & -3.4$\pm$1.8 & -6.0$\pm$6.4  & -8.1$\pm$3.8 & -9.6$\pm$8.1 \\
		14 & -9.5$\pm$3.1 & -17.6$\pm$7.6 & -2.5$\pm$1.4 & -11.9$\pm$8.4 & -7.4$\pm$3.1 & -13.4$\pm$8.1 \\
		15 & -12.8$\pm$5.3 & -13.9$\pm$5.9 & -3.6$\pm$2.1 & -5.2$\pm$4.7 & -10.1$\pm$5.8 & -15.6$\pm$16.0 \\
		16 & -10.0$\pm$3.9 & -20.5$\pm$10.1 & -2.7$\pm$1.7 & -15.0$\pm$10.2 & -8.2$\pm$3.2 & -14.4$\pm$8.3 \\
		17 & -13.1$\pm$5.2 & -16.0$\pm$6.1 & -3.8$\pm$2.2 & -7.3$\pm$5.4 & -10.7$\pm$6.3 & -14.8$\pm$15.8 \\
		18 & -10.5$\pm$4.2 & -26.6$\pm$11.4 & -2.8$\pm$1.7 & -21.7$\pm$11.8 & -8.6$\pm$3.9 & -17.2$\pm$10.7 \\
		19 & -10.1$\pm$3.2 & -11.7$\pm$4.8 & -3.2$\pm$1.7 & -5.6$\pm$5.4 & -7.8$\pm$3.5 & -6.4$\pm$3.7 \\
		20 & -9.4$\pm$3.4 & -15.0$\pm$9.5 & -2.7$\pm$1.5 & -9.2$\pm$10.1 & -7.2$\pm$3.3 & -9.9$\pm$7.5 \\
		21 & -12.5$\pm$5.4 & -14.9$\pm$6.2 & -4.2$\pm$2.3 & -6.8$\pm$8.5 & -10.1$\pm$5.6 & -8.6$\pm$6.8 \\
		22 & -10.3$\pm$3.9 & -17.5$\pm$10.5 & -2.8$\pm$1.7 & -11.3$\pm$11.7 & -8.3$\pm$3.7 & -11.4$\pm$9.0 \\
		23 & -12.5$\pm$5.3 & -13.0$\pm$8.1 & -3.8$\pm$2.1 & -4.5$\pm$9.2 & -9.9$\pm$5.6 & -8.9$\pm$6.1 \\
		24 & -9.7$\pm$3.7 & -18.3$\pm$15.1 & -2.8$\pm$1.5 & -13.4$\pm$14.4 & -7.7$\pm$3.3 & -3.3$\pm$34.3 \\ \hline
		Mean & -10.9$\pm$0.9 & -16.5$\pm$3.1 & -3.2$\pm$0.3 & -9.8$\pm$3.1 & -8.7$\pm$1.2 & -11.3$\pm$8.3 \\ 
		\hline
	\end{tabular}
\end{table}

As expected, the gain of Strategy (1+2+3) compared to Strategy 1-2-3 is, in general, higher than the gain of Strategy (1+2+3) compared to Strategies (1+2)-3 or 1-(2+3), as the latter already work with a type of integration between phases. In addition, the gains obtained by classes 13--23 are generally higher than for classes 1--12, which means that opportunities for improvement are increased by working with a more significant number of product types. 

\subsubsection{Computational performance} \label{sec:desempenho_computacional}

\begin{figure}[!h]
	\centering
	\subfigure[Number of columns of the model, per phase and per set of classes.]{
		{\includegraphics[width=.45\textwidth]{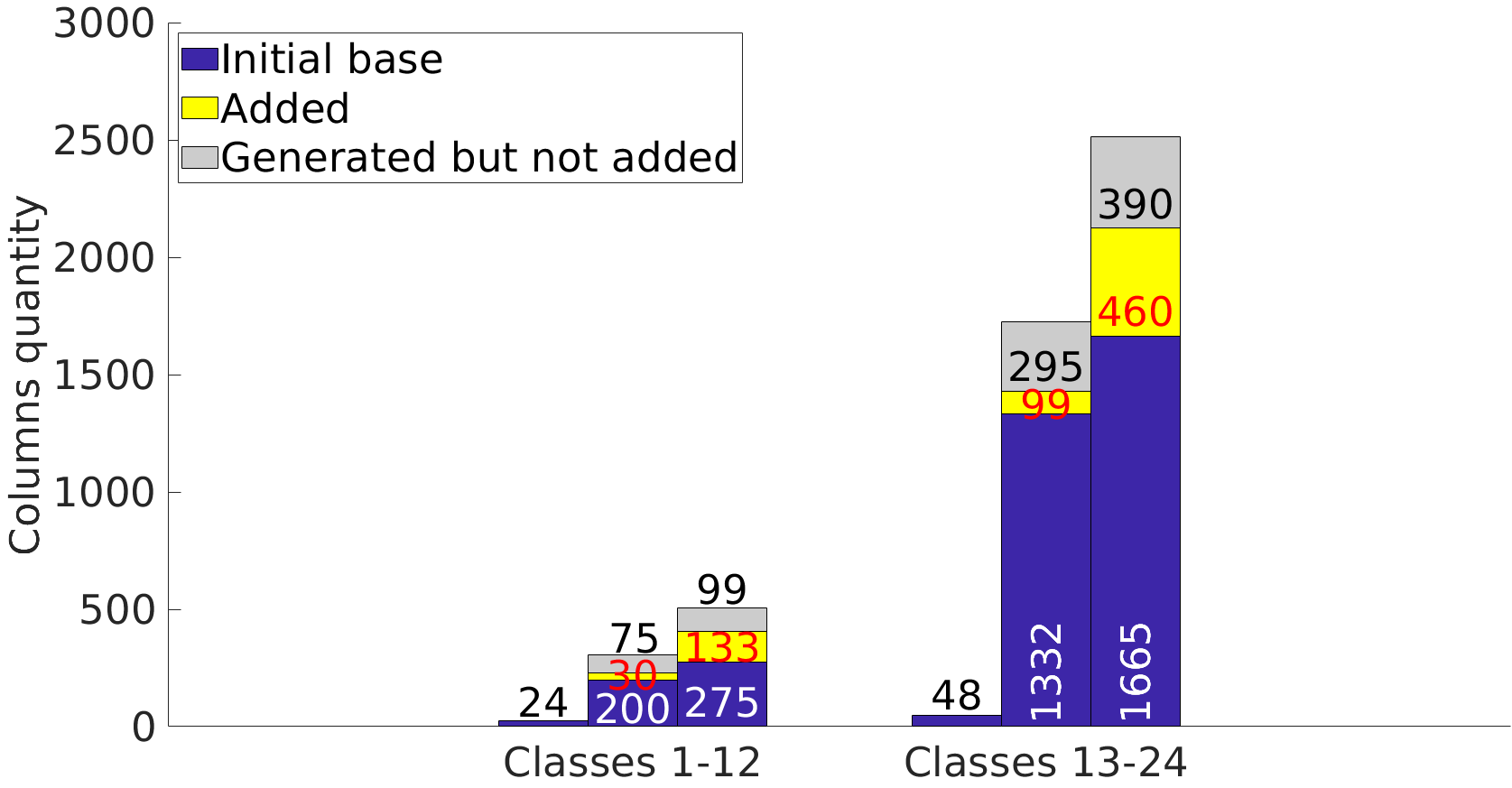}}\label{fig:colunas_classes}}
	\subfigure[Processing time and number of iterations of the heuristic, by class set.]{
		{\includegraphics[width=.45\textwidth]{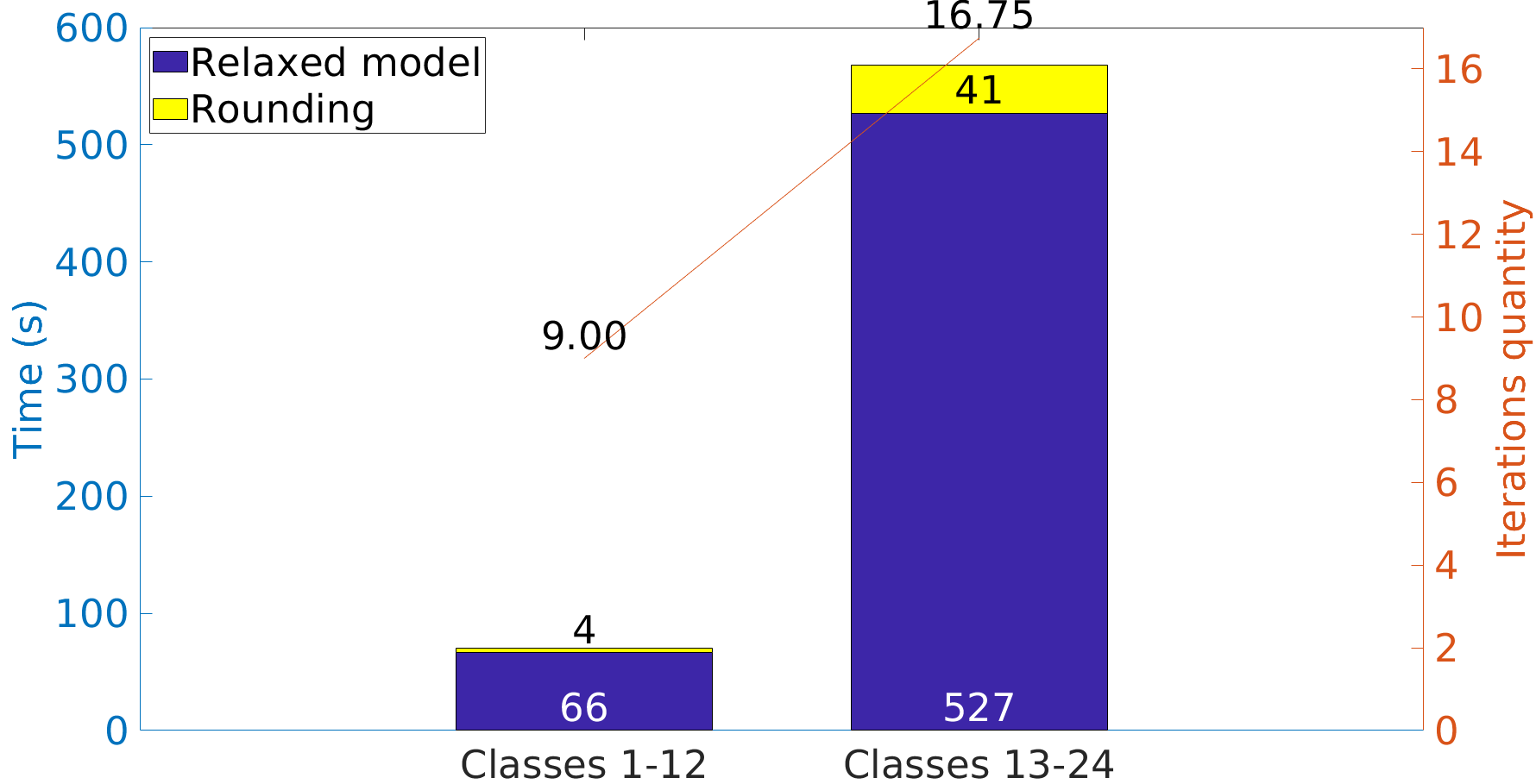}}\label{fig:tempos_classes}}
	\caption{Computational performance.} \label{fig:comp_performance}
\end{figure}

The computational performance of the solution heuristic method was analysed for two sets of classes, those of numbers 1--12, which include three types of jumbos, five types of reels and five types of sheets, and of 13--24, with six types of jumbos, nine types of reels and nine types of sheets. Figure \ref{fig:comp_performance} shows the performance according to different criteria.

Figure \ref{fig:colunas_classes} compares the number of columns of the model per phase for each set of classes. By doubling the number of jumbo types, the number of columns in Phase 1 also doubles, as expected. For Phases 2 and 3, however, by increasing the number of product types from 5 to 9 (ratio of 1.8), the observed increase in the number of columns is more significant---more than 6 times---for the initial base columns, and in the range of 3--4 times for both \textit{added} and \textit{generated but not added} columns. This finding shows the complexity introduced to the model for each type of reel or sheet added. Furthermore, it justifies the use of a specific technique to round the solution obtained by the heuristic method, as it is not possible to solve such large-dimension integer models at once in reasonable time.

In Figure \ref{fig:tempos_classes}, the processing time (primary axis) and the median number of iterations (secondary axis) of the heuristic are compared for each set of classes. The number of iterations increased by 1.86 times, close to the proportion of the increase in the quantity of each type of product, passing from the set of classes 1--12 to 13--24, which suggests a linear relationship between the number of product types and the number of heuristics iterations. For the processing time, however, an increase of 7.98 times in the solution time of the relaxed model and of 10.25 times in the rounding step was verified. This was caused by the more significant number of columns generated and inserted for classes 13--24. 

\subsubsection{Stock cost} \label{sec:analise_estoque_custo}

\begin{figure}[!h]
	\centering
	\subfigure[Phase 2.]{
		{\includegraphics[width=.45\textwidth]{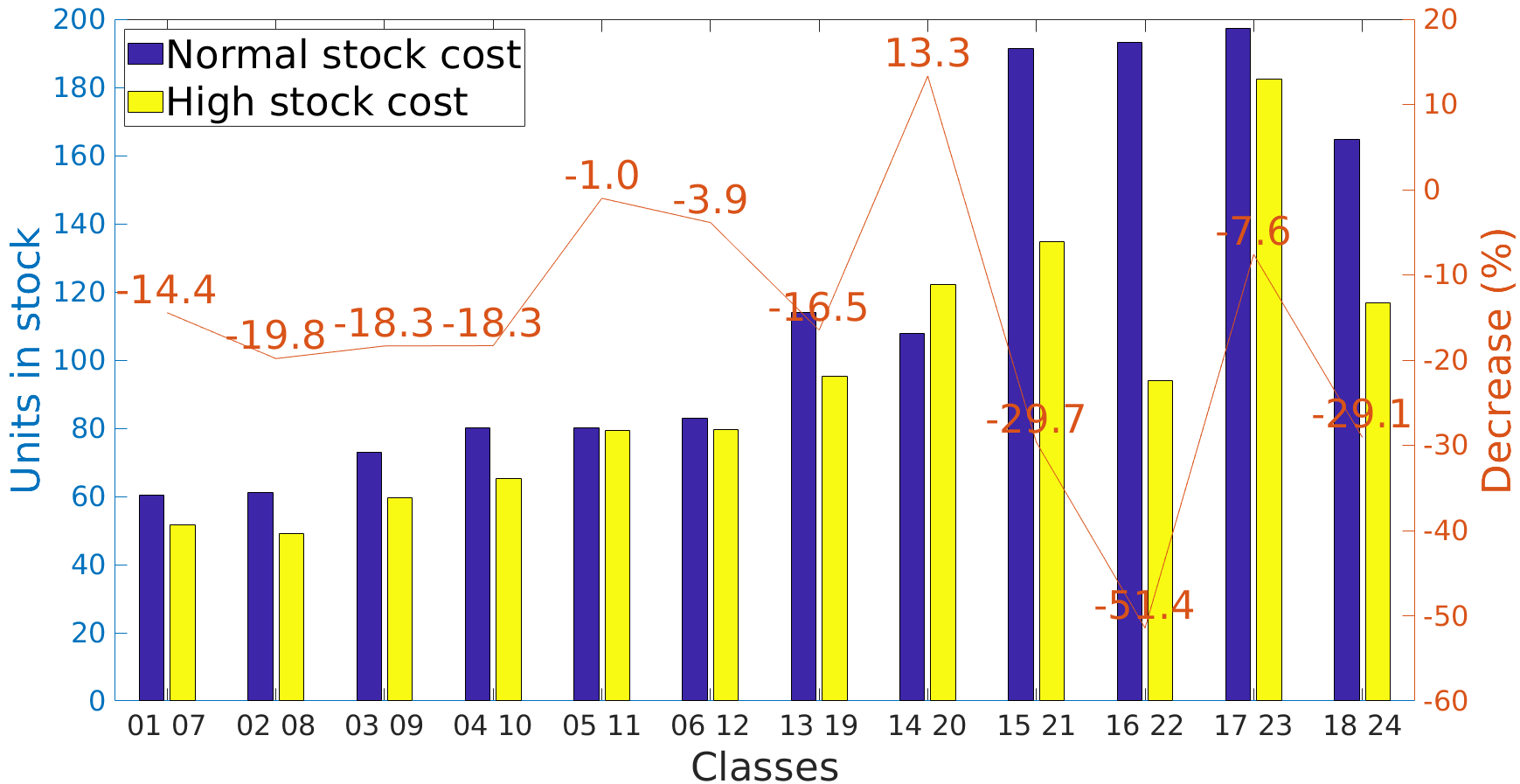}}\label{fig:analise_estoque2}}
	\subfigure[Phase 3.]{
		{\includegraphics[width=.45\textwidth]{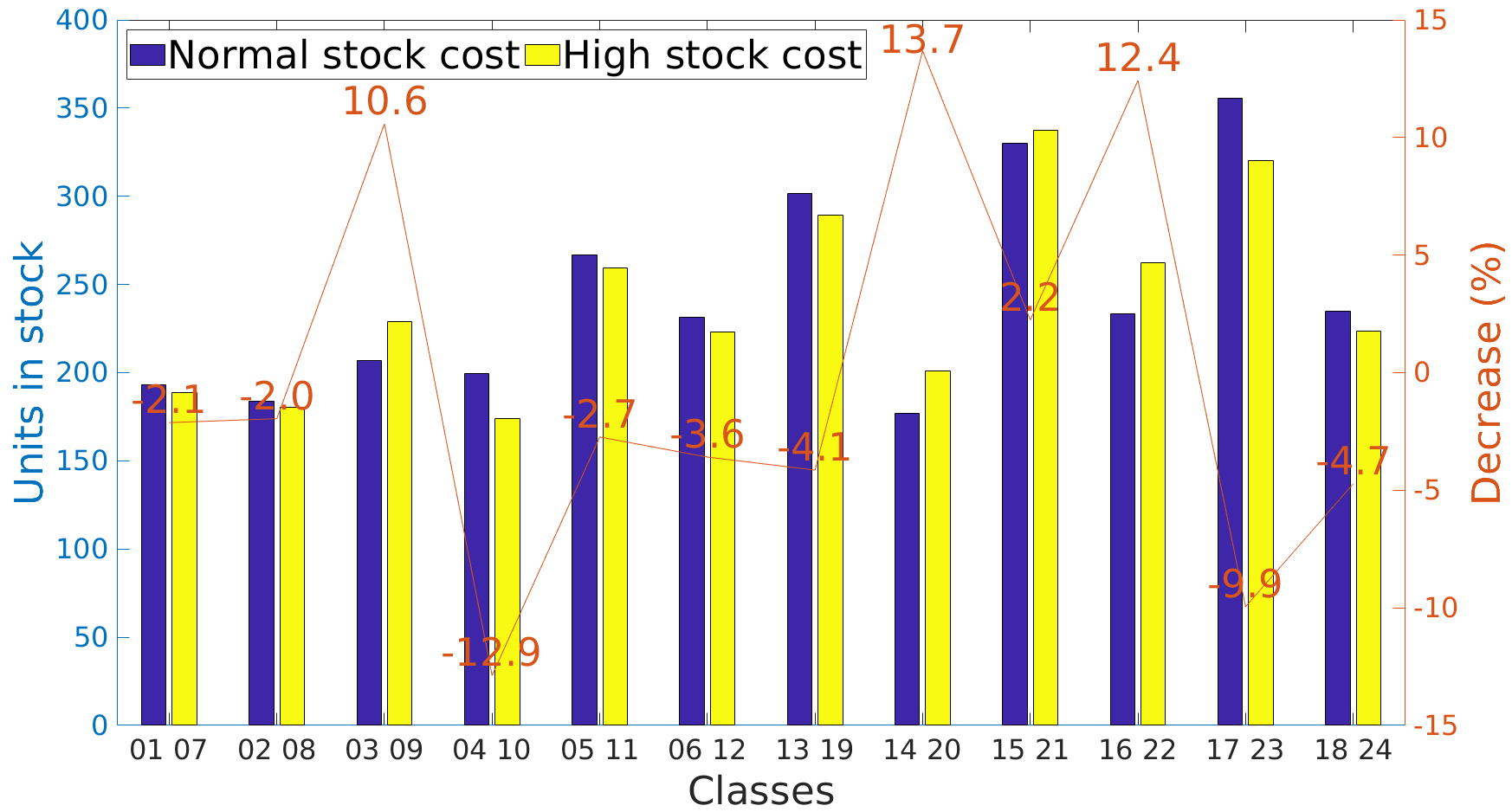}}\label{fig:analise_estoque3}}
	\caption{Number of units in stock for Phases 2 and 3 for normal and high stock cost.} \label{fig:analise_estoque}
\end{figure}

In Figure \ref{fig:analise_estoque}, each pair of bars, corresponding to the median number of units in stock for each phase, has the same settings, except for the stock cost, which varies from normal to high. As expected, there is a drop in the number of units in stock when the cost of inventory increases. Some exceptions do happen for pairs 14--20, 03--09 and 16--22. The reason is that, in the mathematical model, it is not the number of units that is minimised but the weight of these items, which have different storage costs, depending on the type.

\subsubsection{Cutting type} \label{sec:analise_apara}

\begin{figure}[!h]
	\centering
	\subfigure[Material waste.]{
		{\includegraphics[width=.45\textwidth]{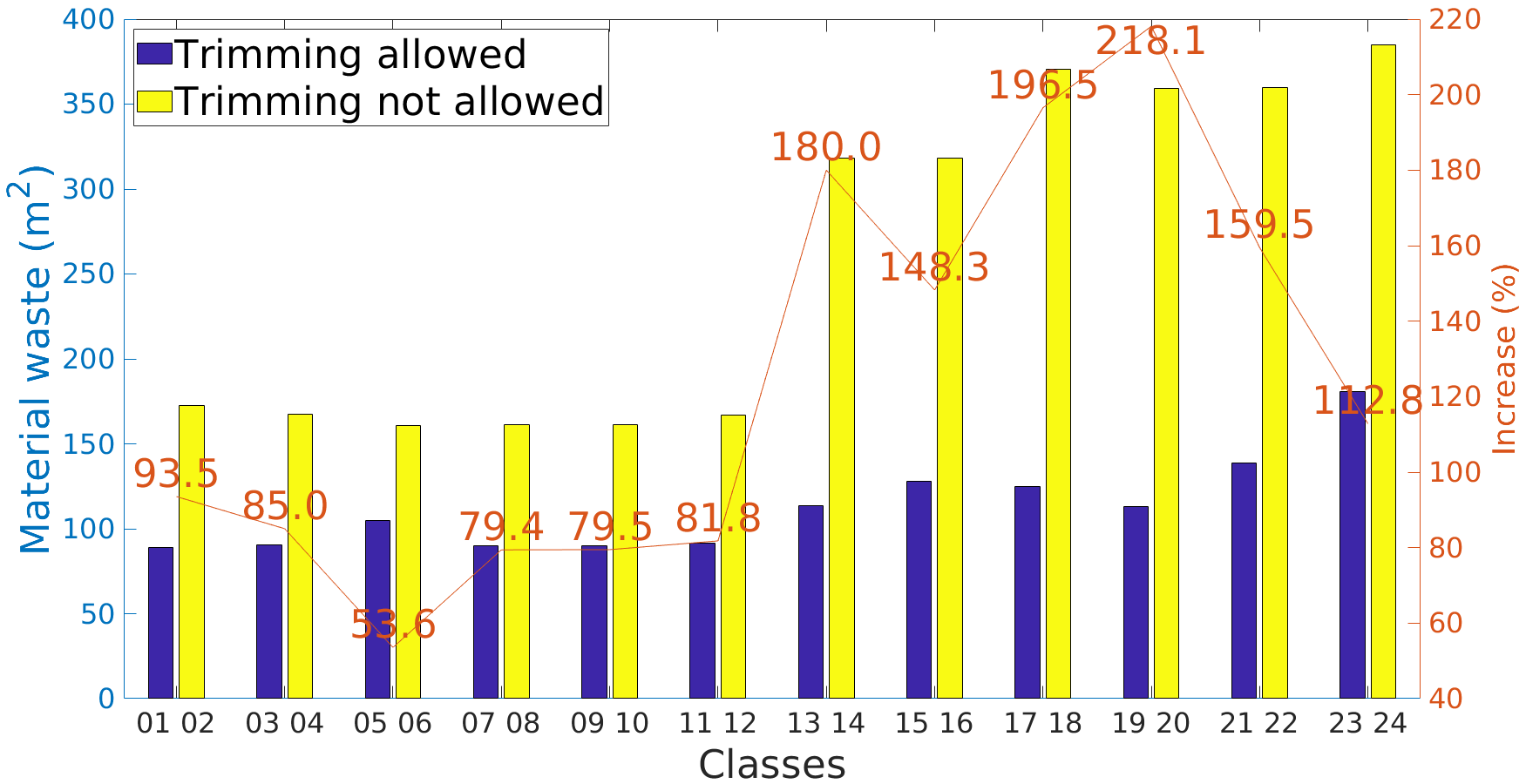}}\label{fig:desperdicio_apara}}
	\subfigure[Number of units in stock.]{
		{\includegraphics[width=.45\textwidth]{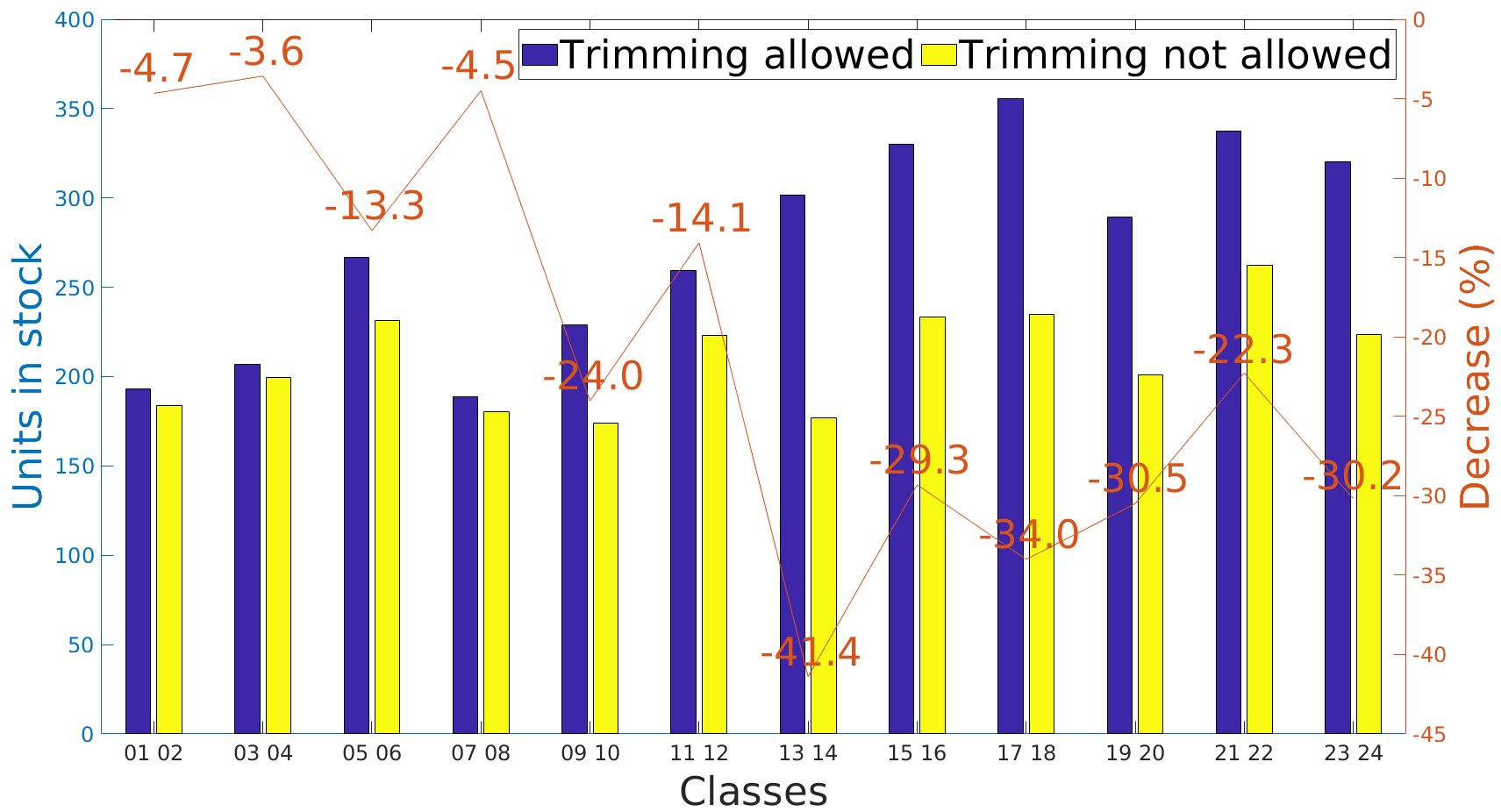}}\label{fig:estoque_apara}}
	\caption{Performance of Phase 3 for the classes with and without trimming.} \label{fig:trimming}
\end{figure}

The results obtained by the odd classes in Table \ref{tab:classes}, which allow trimming in the two-dimensional cutting of reels into sheets, are compared in Figure \ref{fig:trimming} with the even classes, which do not allow this. The material waste in Phase 3 is analysed in Figure \ref{fig:desperdicio_apara}, in which each pair of bars has the same configuration, except for the possibility of two-dimensional cutting with trimmings. For all couples, the waste of material in Phase 3 was smaller with trimming, with more significant reductions for classes 13--24. By allowing cutting with trimming, the algorithm takes more advantage of the object shape, thereby producing less waste. By increasing the number of item types, however, the likelihood of having an item of adequate dimensions to fit the missing space of the pattern also increases, thus reducing waste.

In Figure \ref{fig:estoque_apara}, the number of units of sheets in stock for each pair of bars is compared. Cutting when trimming is not allowed generates less units in stock, with more significant decreases for classes 13--24---with an average 31.3\% of reduction--- than for classes 1--12---with an average 10.7\% of reduction. By decreasing the waste of material through good cutting patterns when trimming is allowed, the model tends to anticipate items that can be cut in the remaining space, thus increasing the number of units in stock.

\section{Conclusions and future works} \label{sec:conclusions}

This article proposed a mathematical model of integer linear optimisation, called the Bi-Integrated Model, that integrates the three phases at the core of the papermaking process. Its objective is to minimise the total costs of production, machine setup and item stocks for each phase. A heuristic method that uses the Simplex Method with column generation for each cutting phase was developed, starting from the relaxed model and rounding the obtained solution to integer values using Relax-and-Fix method.

In analysing the behaviour of costs per phase for each of the four solution strategies, one may note that what makes the Bi-Integrated Model obtain better solutions is not to look for the optimum individually for each phase but rather to seek a better solution for the process as a whole. Thus, sometimes the individual costs of the Bi-Integrated Model were higher than those obtained by the other strategies. However, such solutions produced a more profitable interface with the other phases of the problem, resulting in lower total costs.

The Bi-Integrated Model uses more efficient cutting patterns, resulting in reduced material waste. This is done by anticipating items from subsequent (sub-) periods and inserting more cutting patterns into the master problem compared to the other strategies. However, concerning the number of columns used by the solution, relatively constant values among the strategies were observed. This finding is remarkable because the same costs and machine setup times are maintained. In addition, a small increase in the processing time of the Bi-Integrated Model was verified, which is the price to be paid to support the impressive gain in performance.

For classes with five types of reels and sheets, Strategy (1+2+3) obtained average gains of 11.66\% compared to Strategy 1-2-3, and 8.75\% and 7.26\% compared to Strategies (1+2)-3 and 1-(2+3), respectively. For the classes with nine types of reels and formats, the gains were even more significant: 16.45\%, 9.8\% and 11.25\%, respectively. This finding indicates exciting possibilities for using the Bi-Integrated Model in industrial environments.

Finally, by increasing the cost of inventory, the number of units in stock decreased, as might be expected. In addition, the material waste in Phase 3 is lower when cutting with trimming is allowed, as there is a better use of items in the cutting patterns. However, the number of units in stock is higher when trimming is allowed, as the model anticipates pieces from the following periods. 

Despite our effort to create a generalist model that represents different realities of the paper industry, some simplifications have been made. All paper thicknesses were considered to be unique and neither the diameters of the jumbos nor the diameters of the reels vary. A relevant extension may therefore be to incorporate these variations into the model. 

Regarding the rounding algorithm, although total costs close to the relaxed solution were obtained, there was a restriction on the number of product types manufactured. A \textit{solver} was used to solve mixed-integer problems for each (sub-) period, leading to high processing times. Therefore, an idea for future work may be to develop another rounding technique or to improve the one proposed in this work.

\section*{Funding}
This study was financed in part by FAEPEX-PAPDIC-Unicamp (master’s program scholarships, grants 1123/15 and 519.292-262/15). The second author wishes to thank the Coordenação de Aperfeiçoamento de Pessoal de Nível Superior - Brasil (CAPES) - Finance Code 01P-3717/2017. The third author is grateful to FAEPEX - Unicamp (grants 519.292-285/15 and 519.292-885/15). The fourth author is grateful to FAPESP (grant 2015/02184-7) and FAEPEX-Unicamp (grant 519.292-262/15).

\bibliographystyle{apacite}
\bibliography{manuscript}

\end{document}